\documentclass[11pt]{article}
\usepackage[top=1.1in, left=1in, right=1in, bottom=1.1in]{geometry}
\usepackage{amsmath,amssymb,amsthm,amsfonts,verbatim}
\usepackage{enumerate}
\usepackage{microtype}
\usepackage{url}
\usepackage{algpseudocode}
\usepackage{tikz}
\usetikzlibrary{arrows}
\usepackage{amsmath,amssymb}
\usepackage[justification=centering]{caption}
\usepackage[justification=centering]{subcaption}
\usepackage[english]{babel}
\usepackage{array}
\usepackage{psfrag}
\usepackage{graphicx}
\usepackage{color}
\usepackage{setspace}
\usepackage{hyperref}
\usepackage[utf8]{inputenc}
\usepackage{cleveref}
\usepackage{mathrsfs} 
\usepackage{authblk}
\usepackage{threeparttable}
\usepackage[titletoc,title]{appendix}
\usepackage{xcolor}
\usepackage{algorithm2e}
\usepackage{hyperref}

\newtheorem{remark}{Remark}[section]

\definecolor{red}{rgb}{1,0.2,0.2}

\SetKw{KwInitialization}{Initialization}
\SetKw{KwTheMainBody}{The main body}

\def \mr {\mathbb{R}}
\def \U {\mathcal{U}}

\makeatletter
\newcommand*{\rom}[1]{\expandafter\@slowromancap\romannumeral #1@} 
\makeatother

\hypersetup{colorlinks=true,linkbordercolor=red,linkcolor=blue,pdfborderstyle={/S/U/W 1}}

\def\mr{\mathbb{R}}
\def\mx{\mathbf{x}}
\def\my{\mathbf{y}}
\def\mb{\mathbf{b}}
\def\ma{\mathbf{a}}

\def\U{\mathcal{U}}

\begin{document}

\title{An Ordered Line Integral Method for Computing the Quasi-potential in the case of Variable Anisotropic Diffusion}%

\author[1]{Daisy Dahiya\thanks{ddahiya@math.umd.edu} and Maria Cameron\thanks{cameron@math.umd.edu} }
\affil[1]{Department of Mathematics, University of Maryland, College Park, MD 20742}

\maketitle

\begin{abstract}
Nongradient stochastic differential equations (SDEs) with position-dependent and anisotropic diffusion are often used in biological modeling.
The quasi-potential is a crucial function in the Large Deviation Theory that allows one to estimate transition rates 
between attractors of the corresponding ordinary differential equation and find the maximum likelihood transition paths.
Unfortunately, the quasi-potential can rarely be found analytically. 
It is defined as the solution to a certain action minimization problem.
In this work, the recently introduced Ordered Line Integral Method (OLIM) is extended for computing the quasi-potential for 2D SDEs 
with anisotropic and position-dependent diffusion scaled by a small parameter on a regular rectangular mesh. 
The presented solver employs the dynamical programming principle. 
At each step, a local action minimization problem 
is solved using straight line path segments and the midpoint quadrature rule. 
The solver is tested on two examples where analytic formulas for the quasi-potential are available. 
The dependence of the computational error on the mesh size, the update factor $K$ (a key parameter of OLIMs),
as well as the degree and the orientation of  anisotropy is established.
The effect of anisotropy on the quasi-potential and the maximum likelihood paths is demonstrated on the Maier-Stein model.
The proposed solver is applied to find the quasi-potential and the maximum likelihood transition paths in 
a model of the genetic switch in Lambda Phage between the lysogenic state where the phage reproduces inside the infected cell
without killing it, and the lytic state where the phage destroys the infected cell.
\end{abstract}

{\bf Keywords:}
Quasi-potential, Minimum Action Path, Ordered Line Integral Method, Variable and  Anisotropic Diffusion, Maier-Stein model, Lambda Phage

\section{Introduction}
\label{sec:intro}
Randomness is an integral part of many real-life phenomena such as chemical reactions, biological switches, 
regime changes in climate, and  population dynamics. In order to take it into account, such systems are often modeled 
using Stochastic Differential Equations (SDEs) of the form
\begin{equation}
\label{sde1}
d\mx=\mb(\mx)dt+\sigma(\mx)\sqrt{\epsilon} dW,\quad \mx\in\Omega\subset\mr^d,
\end{equation}
where $\mb(\mx)$ is a twice continuously differentiable vector field, $\epsilon$ is a small parameter, $\sigma(\mx)$ 
is the diffusion matrix, and $W$ is the standard $d$-dimensional Brownian motion. 
In this work, we assume that $\sigma(\mx)$ is invertible everywhere in the domain $\Omega$. 
Along with SDE \eqref{sde1}, we consider the corresponding ODE $\dot{\mx} = \mb(\mx)$.
We assume that the vector field $\mb(\mx)$ is such that 
any trajectory of $\dot{\mx} = \mb(\mx)$ remains in a bounded region and hence approaches some attractor of the system as time goes to infinity.
The noise term $\sigma(\mx)dW$ in SDE \eqref{sde1} enables transitions between the attractors. 

The asymptotic analysis of the dynamics of SDE \eqref{sde1} as $\epsilon\rightarrow 0$ is the subject of the Large Deviation Theory \cite{FW}.
Its key function quantifying the difficulty of the escape from any neighborhood of a given attractor is the quasi-potential.
The quasi-potential can be used to estimate the relative stability of different attractors of $\dot{\mx} = \mb(\mx)$,
the expected transition times between them up to the exponential order \cite{FW}, or, in some cases, sharply \cite{Bouchet}.
The Minimum Action Paths (MAPs)   (a.k.a. Maximum Likelihood transition Paths (MLPs) or instantons) can be readily computed once the quasi-potential is found.
Finally, the invariant probability measure near a given attractor can be approximated using the quasi-potential up to the exponential order.
 
Unfortunately, the quasi-potential  is the solution to
 an action minimization problem that cannot be solved analytically in general. 
 It can be easily shown \cite{quasi} that the quasi-potential $U_E(\mx)$ with respect to a given attractor $E$
 is Lipschitz continuous in any bounded domain but not necessarily differentiable. 
 A well-known example where it is non-differentiable
 is the Maier-Stein model \cite{MS1993,MS1996}. A special case where the quasi-potential is known analytically
 is a linear SDE $d\mx = A\mx dt + \Sigma\sqrt{\epsilon}dW$ where all eigenvalues of the matrix $A$
are negative, and $\Sigma$ is a constant nonsingular matrix \cite{chen,chen1}. 

 The goal of this work is to develop a solver 
for computing the quasi-potential on a regular rectangular mesh for SDE \eqref{sde1} with position-dependent and anisotropic diffusion.
 We have achieved this goal in 2D by extending the recently introduced OLIM-Midpoint 
 algorithm \cite{OLIMs} for uniform isotropic diffusion. 
Our C code \verb|OLIM4VAD.c| containing the presented solver  is posted on M. Cameron's website \cite{Mweb}. 

This work is inspired by our desire to apply the quasi-potential analysis to interesting biological systems.
The particular motivating example is a model for the genetic switch in Lambda Phage from the lysogenic state  where the virus reproduces within an 
infected bacterium Escherichia coli without killing it, to the lytic state where the bacterium is destroyed by the phage. 
The microscopic level model for this switch developed in seminal works by Ackers, Shea, and Johnson in 1980s \cite{AckersJohnShea,SheaAckers} 
 was eventually reduced to
a 2D SDE of the form \eqref{sde1} with a position-dependent diagonal diffusion matrix by Aurell and collaborators  \cite{Aurell,AurellBrown}. 

In the case where the diffusion matrix $\sigma$ in Eq. \eqref{sde1} is constant and symmetric positive definite, 
we have conducted a detailed study of the dependence of the numerical error on the ratio of the eigenvalues of $\sigma$
and the direction of its eigenvectors.  We also devised a test example with a position-dependent matrix $\sigma(\mx)$ 
where the quasi-potential can be found analytically and investigated the dependence of the numerical error on the mesh size $N$ (the mesh is $N\times N$)
and the update factor $K$, a parameter of crucial importance (see \cite{oum1,oum3,quasi,OLIMs} for details).

We have used our solver to compute the quasi-potential and the MAPs for the Maier-Stein model \cite{MS1993,MS1996}
with a family of constant symmetric positive definite diffusion matrices $\sigma$ with eigenvalues $1$ and $2$ and various orientations of its eigenvectors.
 The orientation of the eigenvectors causes a drastic change in the quasi-potential and the MAPs. 
 
 Finally, we applied our solver to compute the quasi-potential and the MAP for the genetic switch model in Lambda Phage \cite{Aurell}.
 Furthermore, we used the Bouchet-Reygner formula to find a sharp estimate for the expected time of the genetic switch 
 from the lysogenic to the lytic state and demonstrated that the variable diffusion matrix proposed in \cite{Aurell} changes the transition rate in comparison with the one computed for the uniform isotropic diffusion. Our results can be compared with earlier ones 
 obtained by other techniques in \cite{zhu,Wang}.
 
 The rest of the paper is organized as follows. In Section \ref{sec:background}, a necessary background on 
 the quasi-potential for SDEs of the form \eqref{sde1} and on methods for computing the quasi-potential on a mesh is provided.
 The quasi-potential solver for SDEs of the form \eqref{sde1} is introduced in Section \ref{sec:method}.
 Numerical tests are presented in Section \ref{sec:test}. The effects of the anisotropic diffusion are illustrated on the Maier-Stein models in Section \ref{sec:MS}.
 An application to the genetic toggle model for Lambda Phage is discussed in Section \ref{sec:GT}.
We provide some additional details for Sections \ref{sec:background} and \ref{sec:GT} in Appendices.


\section{Background} 
\label{sec:background}
In this Section, we provide a necessary background on the quasi-potential for SDE \eqref{sde1},
and explain the significance and challenges of the computation of the quasi-potential on a mesh.

\subsection{Definition and significance of the quasi-potential}
The Freidlin-Wentzell action functional for  SDE \eqref{sde1} is given by
\begin{equation}
\label{FWA}
S_T(\phi) = \frac{1}{2}\int_0^T\|\dot{\phi} - \mb(\phi)\|^2_{A(\phi)}dt,
\end{equation}
where $\phi(t)$ is an absolutely continuous path, $T$ is the travel time, and $A$ is the symmetric positive definite 
matrix defined by $A(\phi(t))=(\sigma(\phi(t))\sigma(\phi(t))^\top)^{-1}$ (note that $\sigma$ is nonsingular everywhere in $\Omega$ by assumption).
Here, $\|\cdot\|_{A(\phi)}$ denotes the norm associated with the matrix $A(\phi)$. 
The quasi-potential at a point $\mx$  with respect to an attractor $E$ of the deterministic system $\dot{\mx}=\mb(\mx)$ is 
the solution of the following action minimization problem:
\begin{equation}
\label{QPdef}
U_E(\mx) = \inf_{T,\phi}\left\{S_T(\phi)~|~\phi(0)\in E,~\phi(T)=\mx,~\phi~\text{is absolutely continuous}\right\}.
\end{equation}
The minimization with respect to the travel time $T$ can be performed analytically resulting in the 
geometric action \cite{hey2,FW} (see Appendix A for details):
\begin{equation}
\label{GeoAct}
S(\psi) =  \int_0^{L}\left(\|\mb(\psi)\|_{A(\psi)} \|\psi_s\|_{A(\psi)} - \langle\psi_s,\mb(\psi)\rangle_{A(\psi)}\right)ds, 
\end{equation}
where $L$ is the length of the path $\psi.$ Thus, the quasi-potential can be defined in terms of the geometric action as
\begin{equation}
\label{ActionMin}
U_E(\mx) = \inf_{\psi}\left\{S(\psi)~|~\psi(0)\in E,~\psi(T)=\mx\right\}.
\end{equation}
Using the Bellman optimality principle \cite{bellman}, it can be shown that the quasi-potential satisfies the Hamilton-Jacobi equation \cite{Bouchet} (see Appendix B for details):
\begin{equation}
\label{HJE_qp}
\langle \nabla U(\mx) , \nabla U(\mx)\rangle_{A(\mx)^{-1}} + 2 \mb(\mx) \cdot \nabla U(\mx) = 0.
\end{equation}
The characteristics $\psi(\alpha)$ of Eq. \eqref{HJE_qp} satisfy
\begin{equation}
\label{char}
\dot{\psi} = \mb(\psi) + A(\psi)^{-1} \nabla U(\psi).
\end{equation}
 Once the quasi-potential is found, one can obtain a MAP going from the attractor $E$ to a given point $\mx$ 
 by shooting the characteristic backward from $\mx$ to $E$, i.e., by integrating Eq. \eqref{char} with reversed parameter $\alpha$:
 \begin{equation}
 \label{char2}
 \dot{\psi} =  - \frac{ \mb(\psi) + A(\psi)^{-1} \nabla U(\psi) }{ \| \mb(\psi) + A(\psi)^{-1} \nabla U(\psi)\|},\quad \psi(0) = \mx.
 \end{equation}
Here the norm $\|\cdot\|$ is the Euclidean 2-norm. The normalization of the right-hand side in Eq. \eqref{char2} allows one to avoid
the slowdown of the integration near equilibria which are often chosen to be the endpoints of the desired MAP.

Suppose that the quasi-potential $U_E(\mx)$ is found. Then one can estimate the expected exit time $\tau_E$ 
from the basin of attraction of $E$ (denoted by $\mathcal{B}(E)$)
and the invariant probability measure $\mu(\mx)$ within any sublevel set of $U_E\subset \mathcal{B}(E)$ up to the exponential order \cite{FW}: 
\begin{equation}
\label{estimates}
\tau_E \asymp e^{\min_{\mx\in\partial\mathcal{B}(E)}U_E(\mx))/\epsilon},\quad \mu(\mx)\asymp e^{-U_E(\mx)/\epsilon}.
\end{equation}
The symbol ``$\asymp$" denotes the logarithmic equivalence.
For the case where the attractor $E$ is an equilibrium point $\mx_0$,  $\arg\min_{\mx\in\mathcal{B}(E)}U_E(\mx) = \mx^{\ast}$ is a saddle point,
and the quasi-potential $U_E(\mx)$ is twice continuously differentiable near $\mx^{\ast}$, a sharp estimate for the expected exit time 
from $ \mathcal{B}(E)$ was derived by Bouchet and Reygner in 2016 \cite{Bouchet}:
\begin{equation}
\label{TransTime}
T _{\epsilon \rightarrow 0} = \frac{2\pi}{{\lambda}^{*}_{+}} \sqrt{\frac{|\det ~H(\mx^{\ast})|}{\det~ H(\mx_0) }} 
\exp\left(\int_{0}^{L} F(\psi(s)) ds\right)\times \exp\left(\frac{U(\mx^*)}{\epsilon}\right).
\end{equation}
Here, ${\lambda}^{*}_{+}$ is the positive eigenvalue of the Jacobian of the vector field $\mb$
at the saddle point $\mx^{\ast}$, 
$\epsilon$ is a small parameter, $H(\mx_0)$ is the Hessian matrix of the quasi-potential at the equilibrium $\mx_0$,
$H(\mx^{\ast})$ is the Hessian at $\mx^{\ast}$, and 
\begin{equation}
F(\mx) := \nabla\cdot \left(\mb(\mx)+\frac{1}{2}A^{-1}\nabla U(\mx)\right) +  \frac{1}{2}\ma(\mx) \cdot \nabla U ,~~~~
\ma_i(\mx) := \sum_{j=1}^d \partial_j (A^{-1})_{ij}(\mx).
\end{equation}
The integral of $F$ is taken along the MAP. 
Here, unlike \cite{Bouchet}, we parametrize the MAP by its arc length rather than the travel time 
because it is more convenient for 
numerical integration. Eq. \eqref{TransTime} is valid for any dimension of the phase space. 
Estimates for the expected exit times for 2D systems were obtained in earlier
works of Maier and Stein in 1990s \cite{MS1993,MS1997}.

\subsection{Significance of computing the quasi-potential on a mesh}
Why does anyone want to compute the quasi-potential in the whole domain $\Omega$ while one can 
compute the MAP using a path-based method such as { MAM \cite{mam}}, GMAM \cite{hey1,hey2} or 
AMAM \cite{zhou1,zhou2} and then find the quasi-potential along it by integrating the geometric action \eqref{GeoAct} or the Freidlin-Wentzell action
\eqref{FWA} respectively?
While finding MAPs by path-based methods is easier in the sense of programming and can be done in arbitrary dimensions, 
path-based methods have some important limitations. First, they might suffer from slow convergence and get stalled without reaching the desired MAP
e.g. in the case where the MAP performs a lot of spiraling. 
Second, the output of a path-based method is always biased by the initial path.
Hence an irrelevant local minimizer can be found instead of the desired global one. 
Third, the maximum likelihood exit location from $\mathcal{B}(E)$
might be unknown in advance. For example, the exit point for the escape to a large loop in the FitzHugh-Nagumo system lying on
the separatrix trajectory was discovered by means of computing the quasi-potential in the whole space \cite{FitzHugh-Nagumo}.

Computing the quasi-potential in the whole domain $\Omega$ naturally resolves all these difficulties.
The global minima of the geometric action are found automatically by the design of the solver.
Shooting MAPs is a straightforward task once the quasi-potential is available. 

Furthermore, the knowledge
of the quasi-potential in the whole space provides much more information about the system than the MAP.
For example, it gives an estimate for the quasi-invariant probability measure (Eq. \eqref{estimates}).
Furthermore, it allows us to decompose the vector field $\mb(\mx)$ into the potential and the rotational components
\begin{equation}
\label{decomp}
\mb(\mx) = -\frac{1}{2}\nabla U_E(\mx) + \mathbf{l}(\mx),\quad{\rm where}\quad \mathbf{l}(\mx)\cdot \nabla U_E(\mx) = 0.
\end{equation}
Note that while $\mb(\mx)$ is assumed to be smooth, $ \nabla U_E(\mx)$ and   $\mathbf{l}(\mx)$ are not necessarily differentiable.
Finally, the quasi-potential computed in the whole region can be used for comparison of stability of various attractors.
A comprehensive overview on this subject in the context of ecological models is given in \cite{nolting}.

{ We have pointed out a number of important advantages of computing the quasi-potential on the mesh  over
using path-based methods for quantifying the asymptotic behavior of a system evolving according to SDE \eqref{sde1}.
The major shortcoming of this approach is that it is limited to low-dimensional systems (2D and 3D) while the  
path-based techniques work in any dimensions. In particular, they can be used for finding
transition paths in systems evolving according to stochastic partial differential equations (see e.g. \cite{hey1}).
Nonetheless, there are plenty of interesting and important low-dimensional biological and ecological models (see e.g. \cite{nolting,FitzHugh-Nagumo})
where the quasi-potential analysis is instrumental. 
Furthermore, in some systems, the effective dynamics is effectively restricted to a low-dimensional 
manifold, and hence the quasi-potential can be computed on it. Our findings on this subject will be reported elsewhere in the future.
}

\subsection{An overview of methods for computing the quasi-potential on a mesh}
The first proposed method \cite{quasi} for the computation of the quasi-potential for 2D
SDEs of the form 
\begin{equation}
\label{sde2}
d\mx = \mb(\mx)dt +\sqrt{\epsilon}dW
\end{equation}
 was based on the Ordered Upwind Method (OUM) \cite{oum1,oum3}. The OUM is
designed for numerical solution of the Hamilton-Jacobi equation of the form $F(\mx,\tfrac{\nabla u}{\|\nabla u\|})\|u\| = 1$.
It employs the dynamical programming principle. The mesh points are computed approximately in the
increasing order of the value function $u$ at them. 
{ The OUM, as well, as well as the OLIM, belong to the family of label-setting methods. 
A comprehensive overview of them is found e.g. in Ref. \cite{CV}. The labels of the mesh points represent their status with respect to 
the running computation. The labels for the OUM are listed below. The labels for the OLIM are borrowed from the OUM.

\begin{itemize}
\item[0.] {\sf Unknown points:} the points where the solution $u$ has not been computed yet, and none of its nearest neighbors is {\sf Accepted} or {\sf AcceptedFront}.
\item[1.] {\sf Considered:} the points { that} have {\sf AcceptedFront} nearest neighbors. 
Tentative values of $u$ that might change as the algorithm proceeds, are available at them.
\item[2.] {\sf AcceptedFront:} the points at which $u$ has been computed and no longer can be updated, and they have at least one {\sf Considered} nearest neighbor.
These and only these points are used for updating {\sf Considered} points.
\item[3.] {\sf Accepted:} the points at which $u$ has been computed and no longer can be updated, 
and they have only {\sf Accepted} and/or {\sf Accepted Front} nearest neighbors.
\end{itemize}
}

At each step, a point with the smallest tentative value of $u$ becomes {\sf AcceptedFront} or {\sf Accepted}, depending on
whether it has or does not have non-{\sf Accepted} nearest neighbors.  
The OUM is designed under the assumption that the
speed function $F(\mx,\tfrac{\nabla u}{\|\nabla u\|})$ is be bounded away from zero and from above leading to a finite anisotropy ratio 
$\Upsilon: = \max_{\mx,\ma}F(\mx,\ma)/\min_{\mx,\ma}F(\mx,\ma)$. $\Upsilon$ is a key parameter of the OUM. It defines the
update radius $\rho = \Upsilon h$, where $h$ is the mesh step. 
Its significance is the following. At each step, every {\sf Considered} or {\sf Unknown} 
point $\mx$ lying within the radius $\rho$ from the new {\sf {\sf Accepted}Front} point is attempted to be updated by the upwind finite difference scheme. 
This scheme is
set up for 
triangles with one vertex $\mx$ and the other two vertices being all possible {\sf {\sf Accepted}Front} points, nearest neighbors of each other, and lying within the distance $\rho$ from $\mx$. 
It was proven in \cite{oum3}
that then the numerical solution of $F(\mx,\tfrac{\nabla u}{\|\nabla u\|})\|\nabla u\| = 1$ converges to the true first arrival viscosity solution. 
On the other hand, it was demonstrated in \cite{oum1,oum3} that a numerical solution computed 
with a smaller update radius may remain bounded 
but fail to converge to the true solution as $h\rightarrow 0$.

The major difficulty in adjusting the OUM for computing the quasi-potential  
was caused by the fact that the anisotropy ratio for the Hamilton-Jacobi equation $\|\nabla U\|^2 + 2\mb(\mx)\cdot\nabla U = 0$ for SDE \eqref{sde2}
was unbounded. 
From now on, we will assume that the attractor $E$ is fixed { and denote the quasi-potential with respect to $E$ simply by $U$
instead of $U_E$.}
If one recasts  $\|\nabla U\|^2 + 2\mb(\mx)\cdot\nabla U = 0$ to the form $F(\mx,\tfrac{\nabla u}{\|\nabla u\|})\|\nabla u\| = 1$, the speed function 
{ $ [2\mb\cdot\tfrac{\nabla u}{\|\nabla u\|}]^{-1}$ }
reaches infinity wherever the 
computation of the quasi-potential  reaches a trajectory starting at the boundary of $\mathcal{B}(E)$
and approaching another attractor as $t\rightarrow \infty$.  This difficulty was overcome in \cite{quasi} by setting up a finite update radius equal to $\rho = Kh$
where $K$ is a positive and reasonably large integer (e.g. $K=20$) and showing that the additional numerical error due to the
 insufficient update radius decayed quadratically with the mesh refinement for the problem of computing the quasi-potential. 
The OUM-based quasi-potential solver proposed in \cite{quasi} became the core of the publicly available 
 R-package {\tt QPot} \cite{cran,Rjournal}.  
 
 The next technical advancement in computing the quasi-potential  for SDEs of the form \eqref{sde2}
 on a mesh was achieved in \cite{OLIMs}
 where a family of the Ordered Line Integral Methods (OLIMs) was introduced. From the OUM \cite{oum1,oum3}, the OLIMs inherit the 
 general structure. 
 From the OUM-based quasi-potential solver \cite{quasi}, 
they inherit the finite update radius set up by brute force. 
 Contrary to the OUM, the OLIMs completely abandon the consideration of the Hamilton-Jacobi equation. Instead,
they focus on the direct finding of the solution of the action minimization problem.   At every step of OLIMs, a local action minimization problem is 
solved involving an approximation of the geometric action by a second-order (or higher) quadrature rule.
 This innovation led to a tremendous reduction (by two to three orders of magnitude) 
 of the numerical error, mainly the error constant, because it largely eliminated the major contribution to the numerical error due to the integration along rather long line segments.
 Like the OUM, the OLIMs are first-order methods due to the use of linear interpolation. 
 However, they are significantly more accurate and tend to have higher effective convergence rates.
Furthermore, a hierarchical update strategy introduced in OLIMs made them four times faster. 
The overall speed-up in OLIMs compared to the OUM, 
taking into account more expensive solutions of local minimization problems, was approximately by the factor of two.
The graphs of the numerical error versus the  CPU time presented in \cite{OLIMs}  eloquently showed the superiority
of the OLIM quasi-potential solvers over the OUM-based one.

The optimal OLIM in terms of the balance between the numerical error and the CPU time is the OLIM-Midpoint
employing the midpoint quadrature rule. In the present work, we extend the OLIM-Midpoint for 
finding the quasi-potential for SDE \eqref{sde1} with position-dependent and anisotropic diffusion matrix $\sigma(\mx)$.
This extension will facilitate the application of  the quasi-potential analysis to a broader class of biological models.  
 
 All methods discussed in this section have been implemented only in 2D. 
 The extension of OLIMs to 3D is currently under development and will be presented elsewhere in the near future.


\section{The Ordered Line Integral Method for Variable and Anisotropic Diffusion}
\label{sec:method}
In this Section, we present an extension of the OLIM-Midpoint
for computing the quasi-potential on 2D meshes for SDEs  of the form \eqref{sde1}. 
The general framework of OLIMs is described in details in \cite{OLIMs}. Here we will briefly go over it and  
focus on the upgrade of the OLIM-Midpoint for position-dependent and anisotropic diffusion, i.e., for the case where 
$\sigma(\mx)$ in Eq. \eqref{sde1} is an arbitrary nonsingular matrix  function   in $\Omega\subset\mr^2$.
For brevity, we will refer to the solver presented here as the {\tt olim4vad} which stands for the Ordered Line Integral Method for
Variable and Anisotropic Diffusion.

%
\subsection{A brief description of the OLIM general framework}
\label{sec:framework}
The quasi-potential is computed by solving the minimization problem \eqref{ActionMin} 
on a regular rectangular mesh. 
We start the algorithm with an initial set of mesh points where the quasi-potential is known, 
and successively update the neighboring points
until the boundary of the computational domain is reached. If the vector field $\mb(\mx)$  has a small
rotational component,
the computation can be continued throughout the whole domain. 
{ Otherwise, it is safer to terminate the computation as soon as a boundary mesh point becomes {\sf Accepted}.}
The nearest neighborhood of a mesh point $\mx_{i,j}$
consists of eight mesh points surrounding it: 
$$
\mathcal{N}(\mx_{i,j}) =  \{\mx_{k,l}~|~|k-i|\le 1,~|l-j|\le 1,~|k-i|+|l-j|>0\}.
$$
Its far neighborhood consists of all mesh points lying within the update radius from it, i.e., within the distance $Kh$ where $K$ is the update factor  (the user-supplied integer)
and $h = \max\{h_1,h_2\}$ where $h_1$ and $h_2$ are the mesh steps in the directions $x_1$ and $x_2$ respectively.
Every mesh point has a state variable attached to it that takes one of four values at each moment of time: {\sf Unknown}, {\sf Considered}, {\sf {\sf Accepted}Front}, or {\sf Accepted}.
Originally, all mesh points are {\sf Unknown}, and the quasi-potential value at them is set to infinity. 
During the initialization, 
the status of the initialized mesh points changes to {\sf Considered}. Then the while-cycle of the OLIM starts. Each step of the while-cycle consists of the following three sub-steps.
First, a {\sf Considered} point $\mx_0$ with the smallest tentative value of the quasi-potential becomes {\sf AcceptedFront}, and
every {\sf AcceptedFront} point that no longer has {\sf Considered} or {\sf Unknown} nearest neighbors becomes {\sf Accepted}. 
Second, 
for each {\sf Considered} point $\mx$ lying within the update radius from the new {\sf {\sf Accepted}Front} point $\mx_0$, proposed update 
values are computed using the one-point update and the triangle updates involving $\mx_0$ (see Section \ref{sec:update}). 
If any of them is smaller than the current tentative value $U(\mx)$,
it replaces it. Third, each nearest {\sf Unknown} neighbor $\mathbf{x}$ of $\mx_0$ becomes {\sf Considered}. One-point update values are computed at $\mx$ using all
{\sf {\sf Accepted}Front} points within the update radius from $\mathbf{x}$. Suppose the minimal of them comes from an {\sf {\sf Accepted}Front} point $\mathbf{y}$.
Then triangle updates at $\mathbf{x}$ are computed from all triangles $(\mathbf{y},\mathbf{z},\mathbf{x})$ where $\mathbf{z}$ is an {\sf {\sf Accepted}Front} 
nearest neighbor of $\mathbf{y}$. This update strategy for computing tentative values at the points changing their  status from {\sf Unknown} to {\sf Considered}
 is named the hierarchical update strategy in \cite{OLIMs}. It was shown that it might lead to a very small increase of the numerical error 
 (by less than one percent) while makes the CPU time four times smaller.
 
{ The algorithm is summarized in the following pseudocode.

 \begin{algorithm}[H]
\KwInitialization{
Start with all mesh points being {\sf Unknown}. 
Compute tentative values { of $U$} at mesh points near the attractor $E$ and make them {\sf Considered} (see Section \ref{sec:init}).
}

\KwTheMainBody\\
\While { the boundary of the mesh is not reached {\bf and} the set of {\sf Considered} points is not empty}{
 {\bf 1:} Make the  {\sf Considered} point $\mx_0$ with the smallest tentative value of $U$ {\sf AcceptedFront}. \\
{\bf 2:} Make all  {\sf AcceptedFront} nearest neighbors of $\mx_0$ that no longer have {\sf Considered} nearest neighbors  {\sf Accepted}.\\
{\bf 3:} \For{ all  {\sf Considered} points $\mx$ within the distance $Kh$ from $\mx_0$}{
{\bf 3a:} Compute the one-point update values ${\sf Q}_{1pt}(\mx_0,\mx)$
 (see Section \ref{sec:1ptu}).\\
{\bf 3b:} Compute the triangle update values ${\sf Q}_{\Delta}(\mx_0,\mx_1,\mx)$ 
for all {\sf AcceptedFront} nearest neighbors $\mx_1$
of $\mx_0$ (see Section \ref{sec:2ptu}).\\
}
{\bf 4:} \For{ all {\sf Unknown} nearest neighbors $\mx$ of $\mx_0$}{
{\bf 4a:} Change status of $\mx$ to {\sf Considered} \\
{\bf 4b:} \For{ all {\sf Accepted Front} points $\my$ lying within the distance $Kh$ from $\mx$}{
Compute the one-point update values ${\sf Q}_{1pt}(\my,\mx)$  (see Section \ref{sec:1ptu}).
}
{\bf 4c:}  Find the  minimizer $\my_0$ of ${\sf Q}_{1pt}(\my,\mx)$.\\
{\bf 4d:} \For{ for all {\sf AcceptedFront} nearest neighbors $\my_1$ of $\my_0$}{
Compute the triangle update values ${\sf Q}_{\Delta}(\my_0,\my_1,\mx)$
 (see Section \ref{sec:2ptu}).
}
}
}
\caption{The outline of OLIMs with the hierarchical update}
\end{algorithm}
 }

 Next, we will elaborate the initialization and the computation of the updates specific for {\tt olim4vad}.

\subsection{Initialization}
\label{sec:init}
There are three types of attractors in 2D: asymptotically stable equilibria, stable limit cycles, and  sets consisting of equilibrium points and 
heteroclinic or homoclinic trajectories. 
The initialization of a mesh point $\mx$ near the last two types of attractors can be done as follows.
Find the value of the  geometric action
for all straight line segments connecting $\mx$ with the points $\my$ representing the attractor such that $\|\mx - \my\|\le Kh$
and select the minimal of them to be $U(\mx)$.

The initialization near an asymptotically stable equilibrium can be done using the 
analytic formula for the quasi-potential for a linear 2D SDE of the form \eqref{sde2} with a constant diffusion matrix \cite{quasi}. 
In order to do this, one needs to make a variable change reducing SDE \eqref{sde1} 
 to SDE \eqref{sde2}. 
Let us elaborate this procedure. The linearization near the asymptotically stable equilibrium $\mx_0$ 
gives: 
\begin{align}
\mb(\mx) &= J(\mx - \mx_0) + O(\|\mx-\mx_0\|^2),~~{\rm where}~~ 
J = \left(\left.\frac{\partial b_i(\mx)}{\partial x_j}\right\vert_{\mx = \mx_0}\right)_{i,j = 1,2},\label{linb}\\
\sigma(\mx)& = \Sigma + O(\|\mx - \mx_0\|),~~{\rm where}~~ \Sigma: = \sigma(\mx_0). \label{linsig}
\end{align}
Multiplying SDE \eqref{sde1} by $\Sigma^{-1}$, introducing a new variable $\mathbf{y} = \Sigma^{-1}(\mx - \mx_0)$, and 
keeping only the largest order terms in the deterministic and the stochastic parts of the right-hand side we get:
\begin{equation}
\label{lineq}
d\mathbf{y} =[ \Sigma^{-1}J\Sigma]\mathbf{y} dt + \sqrt{\epsilon}dW.
\end{equation}
The quasi-potential decomposition for the linear SDE \eqref{lineq} is given by the formula
 \begin{align}
& \U(\my) =\my^\top \left[\begin{array}{cc} \mathcal{A}&\mathcal{B}\\\mathcal{B}&\mathcal{C}\end{array}\right] \my,
 \quad{\rm where}\label{cam}\\
&\mathcal{A} = -(\alpha g_{11} + \beta g_{21}),\quad
\mathcal{B} = -(\alpha g_{12} + \beta g_{22}),\quad
 \mathcal{C} = -(\alpha g_{22} - \beta g_{12}),
\notag \\
&\alpha = \frac{(g_{11} + g_{22})^2}{(g_{11} + g_{22})^2 + (g_{21} - g_{12})^2},\quad 
\beta = \frac{ (g_{21} - g_{12})(g_{11} + g_{22})}{(g_{11} + g_{22})^2 + (g_{21} - g_{12})^2},\notag
\end{align}
where $g_{ij}$, $i,j=1,2$, are the entries of the matrix $\Sigma^{-1}J\Sigma$.
Returning to the variable $\mx$ we find the quasi-potential near the asymptotically stable  equilibrium $\mx_0$:
\begin{equation}
\label{qx}
U(\mx) \approx (\mx-\mx_0)^\top\Sigma^{-\top} \left[\begin{array}{cc} \mathcal{A}&\mathcal{B}\\
\mathcal{B}&\mathcal{C}\end{array}\right] \Sigma^{-1}(\mx - \mx_0)\equiv (\mx-\mx_0)^\top M (\mx-\mx_0).
\end{equation}
We will refer to the matrix $M$ in Eq. \eqref{qx} as the \emph{quasi-potential matrix}.
Eq. \eqref{qx} is used to initialize four nearest mesh points around $\mx_0$
 if $\mx_0$ is not a mesh point, or eight nearest neighbors of $\mx_0$ otherwise.

\begin{remark}
Note that, slightly counterintuitively, 
the { quasi-potential} matrix  for the linear SDE $d\mx = J\mx + \Sigma\sqrt{\epsilon}dW$ is {\bf not} equal,  in general, 
to $\Sigma^{-\top}(\Sigma^{-1} V\Sigma)\Sigma^{-1} \equiv (\Sigma\Sigma^\top)^{-1}V$ 
where $ V $ is the { quasi-potential} matrix for the linear SDE $d\mx = J\mx + \sqrt{\epsilon}dW$.
Indeed, let $J = -V+L$ where $V$ is symmetric positive definite and $VL$ is antisymmetric.
Then $\Sigma^{-1}J\Sigma = -\Sigma^{-1} V\Sigma+\Sigma^{-1} L\Sigma$, but
$\Sigma^{-1} V\Sigma$ is not symmetric unless $\Sigma$ is a multiple of an orthogonal matrix. 
\end{remark}


\subsection{Update Rules}
\label{sec:update}
The OLIMs involve two types of update rules: one-point update and triangle update. Here we elaborate them for {\tt olim4vad}. 

\subsubsection{One-point update}
\label{sec:1ptu}
Let $\mx_0$ be an {\sf {\sf Accepted}Front} mesh point,  and  $\mx$ be a {\sf Considered} point lying within the update radius $Kh$ from $\mx_0$.
Let $U_0$ and $U(\mx)$ be the values of the quasi-potential at $\mx_0$ and $\mx$ respectively.
The one-point update at $\mx$ from $\mx_0$ is given by
\begin{align}
{\sf Q}_{1pt}(\mx_0,\mx) &= U_0  +\|\mx - \mx_0\|_{A_m}\|\mb(\mx_m)\|_{A_m} - 
\langle\mx-\mx_0,\mb(\mx_m)\rangle_{A_m},\label{1pt}\\
&{\rm where}\quad \mx_m: = \frac{\mx+\mx_0}{2},\quad A_m\equiv A(\mx_m); \quad {\rm then} \notag\\
U(\mx) &= \min\{ {\sf Q}_{1pt}(\mx_0,\mx) ,U(\mx)\}.\label{1ptu}
\end{align}
${\sf Q}_{1pt}(\mx_0,\mx)$ is the approximation of the geometric action \eqref{GeoAct} along the straight line segment $[\mx_0,\mx]$
by the midpoint quadrature rule. Eq. \eqref{1ptu} indicates that the proposed update value 
 ${\sf Q}_{1pt}(\mx_0,\mx)$ replaces the current tentative value at $\mx$
if and only if ${\sf Q}_{1pt}(\mx_0,\mx) < U(\mx)$.

\subsubsection{ Triangle update}
\label{sec:2ptu}
\begin{figure}[htbp]
\begin{center}
\includegraphics[width=0.5\textwidth]{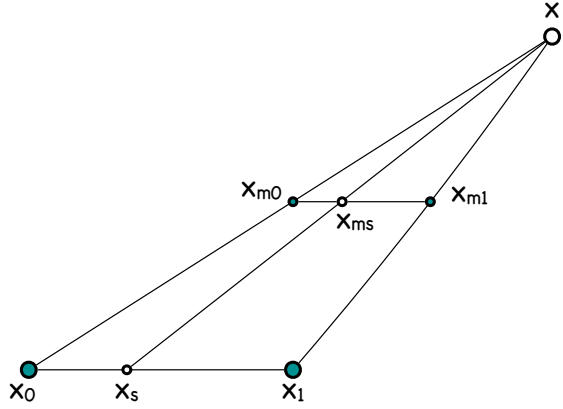}
\caption{An illustration for the triangle update.}
\label{fig:triup}
\end{center}
\end{figure}

Let $\mx_0$ and $\mx_1$ be two {\sf {\sf Accepted}Front} mesh points, nearest neighbors of each other,  with quasi-potential values $U_0$ and $U_1$ respectively.
Let $\mx$ 
be a {\sf Considered} mesh point with a current tentative value $U(\mx)$  lying within the update radius $Kh$ from $\mx_0$ or $\mx_1$ (see Fig. \ref{fig:triup}).
Assume that the points $\mx_0$, $\mx_1$ and $\mx$ form a nondegenerate triangle.
The proposed update value for $U(\mx)$ from the triangle  $(\mx_0,\mx_1,\mx)$ is the solution of the following constrained minimization problem:
\begin{align}
{\sf Q}_{2pt}(\mx_0,\mx_1,\mx)&=\min_{s \in [0,1]} \left\{ U_0 + s(U_1-U_0) \right. \notag \\
&+\left. 
\|\mx-\mx_s\|_{A_{ms}} \|\mb_{ms}\|_{A_{ms}} -\langle\mx-\mx_s,\mb_{ms}\rangle_{A_{ms}} \right\}, \label{2pt_update}\\
U(\mx) &= \min\{ {\sf Q}_{2pt}(\mx_0,\mx_1,\mx) ,U(\mx)\}.\notag
\end{align}
The point $\mx_s:=\mx_0+s(\mx_1 - \mx_0)$ lies on the segment $[\mx_0,\mx_1]$.
In Eq. \eqref{2pt_update}, $A_{ms}$ and $\mb_{ms}$ denote the 
approximations to the matrix $A$ and the vector field $\mb$ respectively at the midpoint of the line segment 
$[\mx_s,\mx]$ found by linear interpolation: 
\begin{align}
\mb_{ms}&=\mb_{m0}+s(\mb_{m1}-\mb_{m0}), ~~{\mbox{and~~}}A_{ms}=A_{m0}+s(A_{m1}-A_{m0}),\label{Abappr}\\
\mb_{m0}&=\mb \left(\frac{\mx_0+\mx}{2} \right),~~~~\mb_{m1}=\mb \left(\frac{\mx_1+\mx}{2} \right),\notag \\
A_{m0}& = A \left(\frac{\mx_0+\mx}{2} \right),~~~~A_{m1}= A \left(\frac{\mx_1+\mx}{2} \right).\notag
\end{align}
 To solve problem \eqref{2pt_update}, we differentiate the function to be minimized with respect to $s$ and
set the derivative to zero. If the signs of this derivative are different at the endpoints $s=0$ and $s=1$, 
we solve the nonlinear 1D equation and find a point $s^{\ast}$, a candidate for the minimizing value of $s\in[0,1]$.
Otherwise, {\sf Q}$_{2pt}$ is not provided from this triangle and $U(\mx)$ remains unchanged.
The nonlinear equation to be solved is
\begin{align}
&U_{1}-U_{0}+ \frac{\|\mb_{ms}\|_{A_{ms}}}{\|\mx - \mx_s \|_{A_{ms}}}\left[ (\mx-\mx_s)^\top A_{ms} \delta\mx
+ \frac{1}{2} (\mx-\mx_s)^\top \delta A_m (\mx-\mx_s) \right]  \notag \\
&+ \frac{\|\mx - \mx_s \|_{A_{ms}}}{\|\mb_{ms}\|_{A_{ms}}} \left[  \mb_{ms}^\top A_{ms}\delta \mb_m + \frac{1}{2} \mb_{ms}^\top \delta A_m \mb_{ms} \right] \notag \\
&- \left( \delta \mx^\top A_{ms} \mb_{ms} + (\mx-\mx_s)^\top A_{ms}\delta\mb_{m}  +(\mx-\mx_s)^\top\delta A_m \mb_{ms}\right) = 0, \label{nonlineq}\\
&{\rm where}~~\delta\mx : = \mx_0-\mx_1,\quad \delta\mb_m : = \mb_{m1}-\mb_{m0},~~\delta A_m: = A_{m1}-A_{m0}.\notag 
\end{align}
We use Wilkinson's hybrid secant/bisection method \cite{wilkinson,stewart} to solve Eq. \eqref{nonlineq}.

The use of the linear interpolation renders {\tt olim4vad} at most first order. 
Since the quasi-potential is not necessarily differentiable, 
and hence its level sets are not necessarily smooth, 
higher order interpolation of $U$ between the {\sf {\sf Accepted}Front} points might lead to larger errors. 
The advantage of the use of the midpoint quadrature rule  is gained due to obtaining  
$O(l^3)$ accurate approximations  for the geometric actions along rather long ($l\sim Kh$) integration paths.
A detailed error analysis would closely repeat the one conducted in \cite{OLIMs}.

\section{Numerical Tests}
\label{sec:test}
In this Section, we test the {\tt olim4vad} on two problems where the quasi-potential is available analytically: 
a linear SDE with a constant non-diagonal diffusion matrix and
a nonlinear SDE with a variable diffusion matrix. 

The linear SDE test problem is designed so that 
the magnitude of the rotational component is approximately 10 times larger than that  of the potential 
component. 
In this example, we study the dependence of the numerical error and 
the optimal choice of the update factor $K$ on the direction of the eigenvectors and the ratio 
of the {eigenvalues} of the diffusion matrix.

In the nonlinear example, the magnitudes of the potential and the rotational components are equal, 
but the diffusion matrix changes  considerably throughout the computational domain.  
For this example, we establish the dependence of the numerical errors on the mesh size and the update factor $K$.

\subsection{A linear SDE with constant anisotropic diffusion}
We consider the SDE
\begin{equation}
\label{linear_test}
d\mx=J \mx dt+\Sigma \sqrt{\epsilon} dW,
\end{equation}
where the matrices $J$ and $\Sigma$ are given by
\begin{equation}
\label{linear_test_input}
J =\left[\begin{array}{cc} -2 & -10  \\ 20 & -1  \end{array}\right], ~~~~ 
\Sigma =\left[\begin{array}{cc} \cos \alpha & -\sin \alpha \\ \sin \alpha & \cos \alpha \end{array}\right] \left[\begin{array}{cc} 1 & 0 \\ 0 & \gamma \end{array}\right]
 \left[\begin{array}{cc} \cos \alpha & \sin \alpha \\ -\sin \alpha & \cos \alpha \end{array}\right].
\end{equation}
The parameters  $\alpha$ and $\gamma$ determine the orientation of the eigenvectors and the ratio of the eigenvalues of $\Sigma$ respectively. 
We consider the following sets of values of $\alpha$ and $\gamma$: 
\begin{equation}
\label{alpha}
\alpha \in \left\lbrace 0, \frac{\pi}{8}, \frac{\pi}{4}, \frac{3\pi}{8},\dots,\frac{7\pi}{8}\right\rbrace
\end{equation}
and 
\begin{equation}
\label{gamma}
\gamma \in \left\lbrace \frac{2}{7}, ~\frac{1}{3},~ \frac{2}{5}, ~\frac{1}{2}, ~\frac{2}{3}, ~1,~ \frac{3}{2},~ 2,~ \frac{5}{2},~ 3,~ \frac{7}{2} \right\rbrace.
\end{equation} 
The exact quasi-potential for SDE \eqref{linear_test} is found by Eqs. \eqref{cam}-\eqref{qx}.

\begin{figure}[htbp]
\begin{center}
\centerline{
(a)\includegraphics[width = 0.4\textwidth]{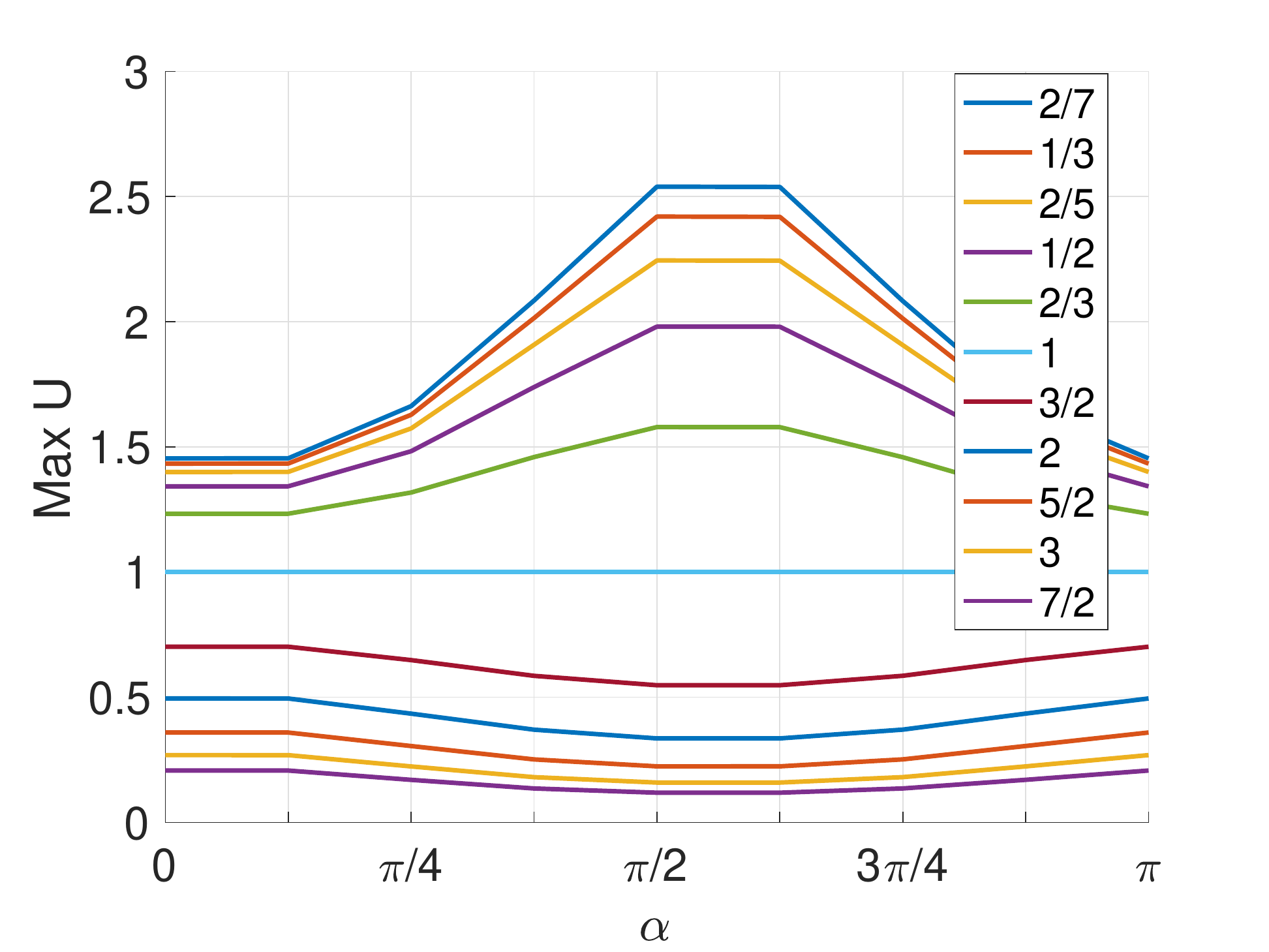}
(b)\includegraphics[width = 0.4\textwidth]{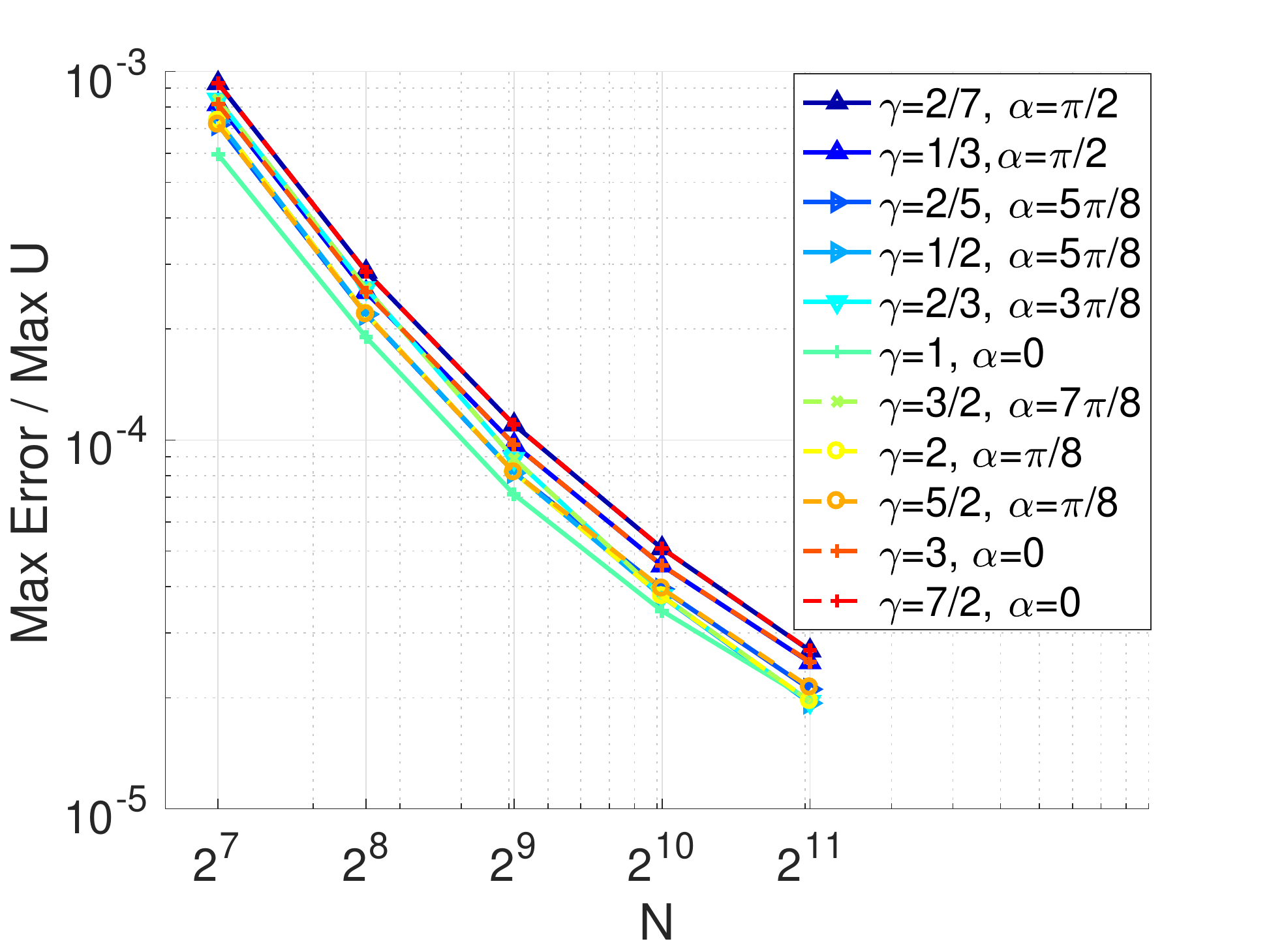}
}
\centerline{
(c)\includegraphics[width = 0.4\textwidth]{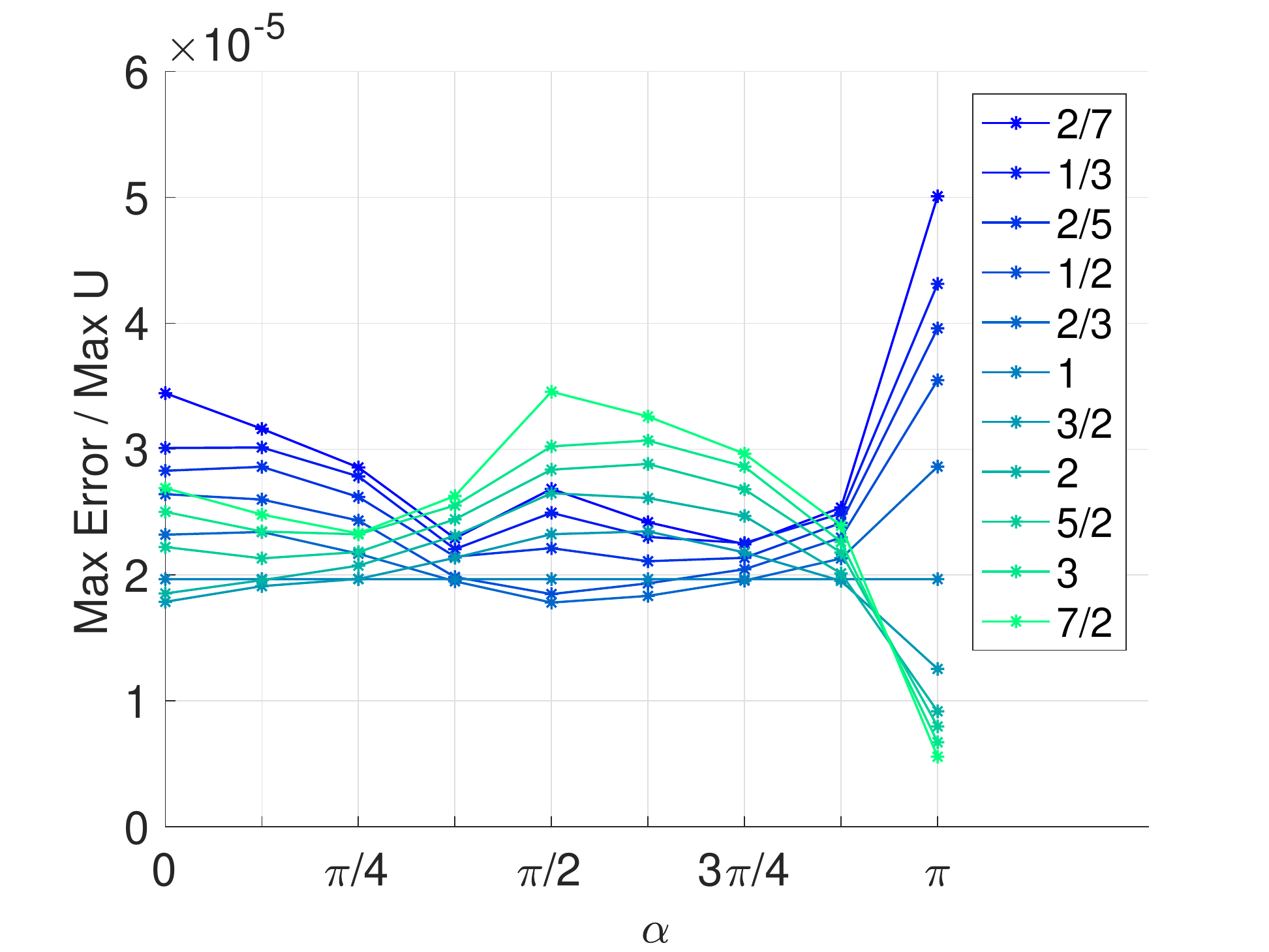}
(d)\includegraphics[width = 0.4\textwidth]{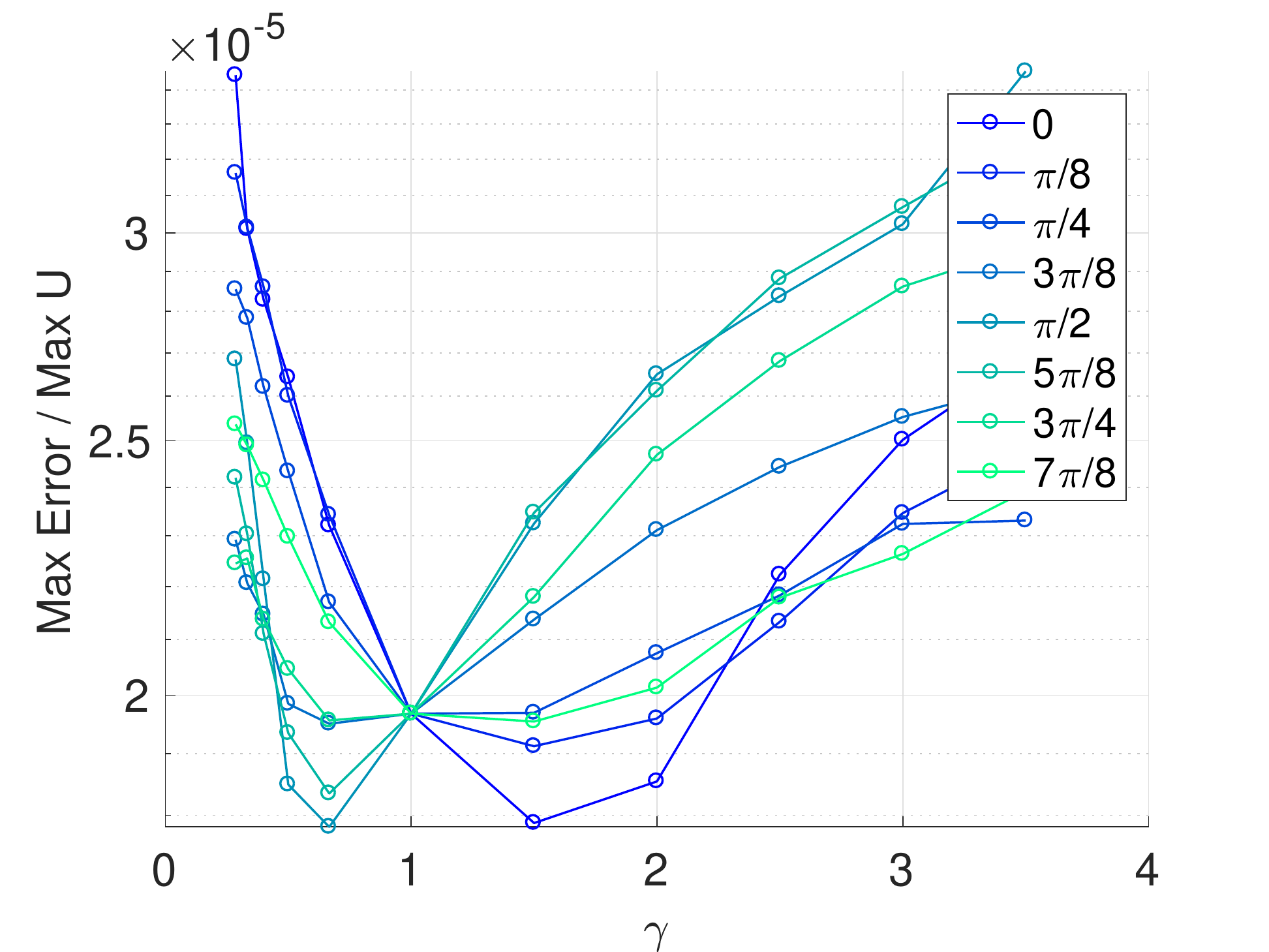}
}
\centerline{
(e)\includegraphics[width = 0.4\textwidth]{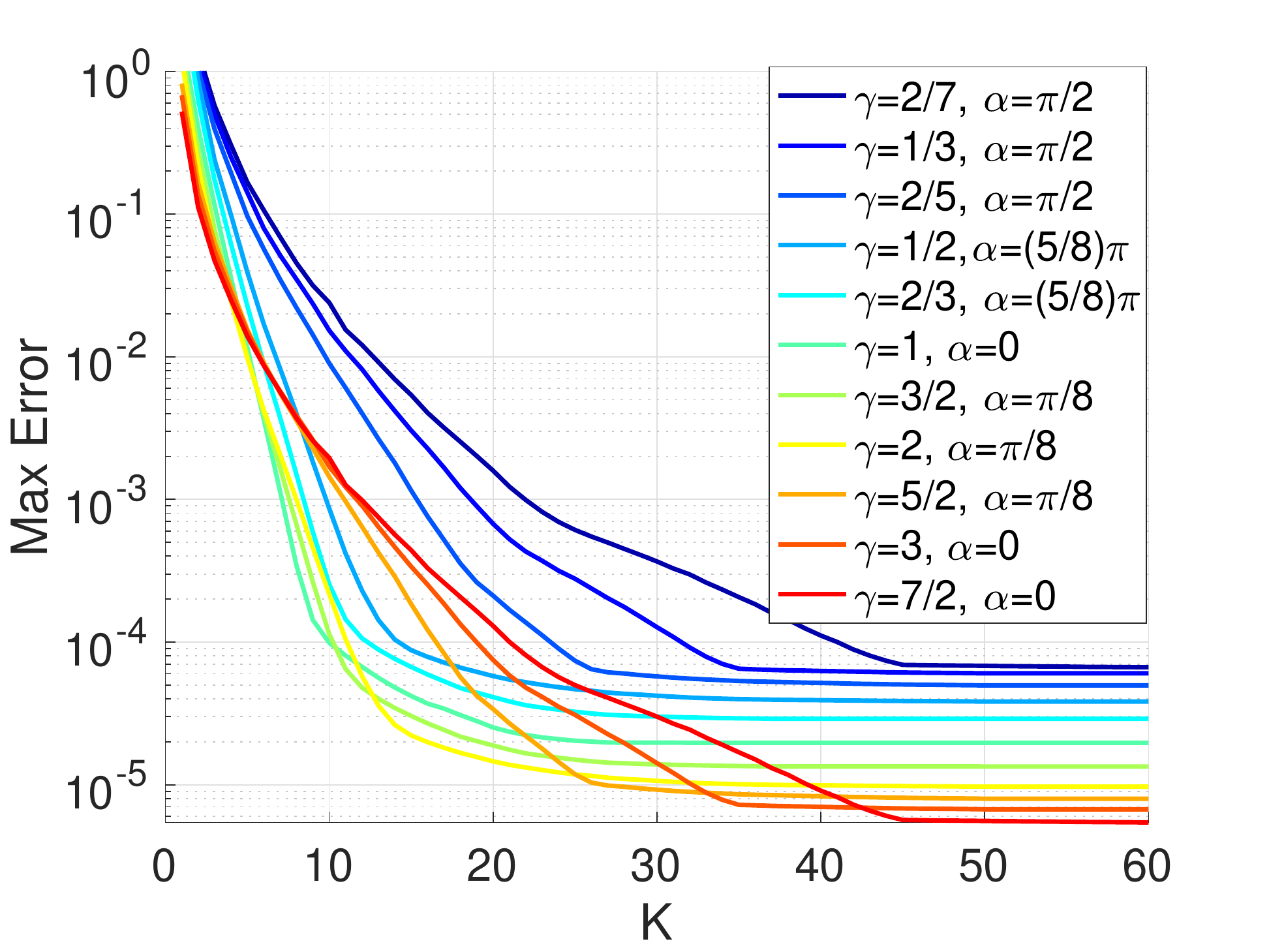}
(f)\includegraphics[width = 0.4\textwidth]{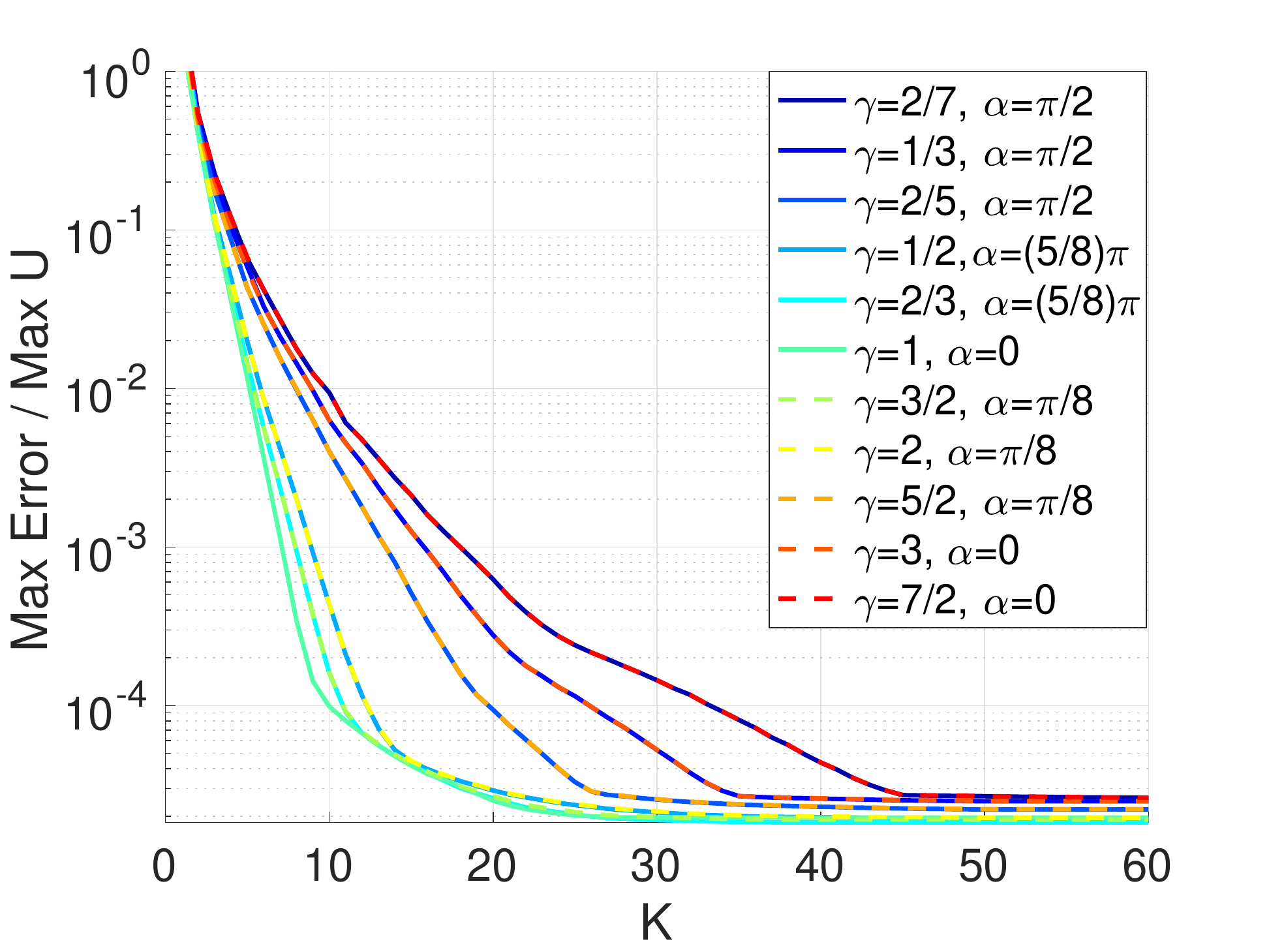}
}
\caption{Numerical errors produced by {\tt olim4vad} for SDE \eqref{linear_test}-\eqref{linear_test_input}.
(a): The maximal value of the computed quasi-potential as the function 
of $\alpha$ for various values of $\gamma$ (see Eq.\eqref{linear_test_input}). 
(b): The normalized maximal absolute error versus $N$ where $N\times N$ is the mesh size.
The graphs are plotted for $K=50$ and for all values of $\gamma$.
(c): The normalized maximum absolute error versus $\alpha$ for each value of $\gamma$. The graphs are given for $N=2048$ and $K=50$.
(d):  The normalized maximum absolute error versus $\gamma$ for each value of $\alpha$. The graphs are given for $N=2048$ and $K=50$.
(e): The maximal absolute error versus the update parameter $K$ for all values of $\gamma$ and some values of $\alpha$. 
(f): The normalized maximal absolute error versus the update parameter $K$ for all values of $\gamma$ and some values of $\alpha$. 
}
\label{fig:linear}
\end{center}
\end{figure}

We have conducted a large number of numerical tests of {\tt olim4vad} on the linear model in order to
$(i)$ ensure the proper decay of the numerical error with the mesh refinement,
$(ii)$ establish the dependence of an optimal choice of the update factor $K$ on the ratio of the eigenvalues of the matrix $\Sigma$ (i.e., on the parameter $\gamma$), and
$(iii)$ establish the dependence of the numerical error on $\gamma$  
and the direction of eigenvectors determined by $\alpha$.
We have measured the maximal absolute errors, the RMS errors, and the CPU times for each combination of values of $\alpha$ and 
$\gamma$ in Eqs. \eqref{alpha} and \eqref{gamma} respectively, and for each combination 
of values of $N = 2^p$, $p = 7,8,9,10,11$ (the mesh size is $N\times N$),
and $K = 1,2,\ldots,60$.  
The computations have been conducted in $[-1,1]^2$ domain starting at the origin, the asymptotically stable equilibrium of $\dot{\mx}=J\mx$,
and terminated as soon as the boundary of the domain was reached. The termination upon reaching the boundary
is necessary due to the fact the rotational component exceeds the potential one by
an order of magnitude making the characteristics spiral quite densely. Hence, they might return back to the square  $[-1,1]^2$ after leaving it.
This would destroy the convergence of the numerical solution to the exact one.

The scale of the quasi-potential is affected by the values of $\gamma$ and $\alpha$.
The maximal values of the quasi-potential $U$ as functions of 
$\alpha$ are plotted for all considered values of $\gamma$ in Fig. \ref{fig:linear}(a).
The maximal value of $U$ reached by our computations ranges from approximately $0.12$ to $2.5$.
Therefore, we divide the maximal absolute errors (i.e., the maximal absolute differences between the computed and the exact solutions)
by the maximal values of the computed quasi-potential and refer to them as the \emph{normalized maximal absolute errors}.

Fig. \ref{fig:linear}(b) displays the dependence of the normalized maximal absolute 
error on $N$ for all values of $\gamma$. The  value of $\alpha$ that maximizes the error is chosen for each $\gamma$.
The least squares fit to $E = CN^{-p}$ gives the values of $p$ ranging from $1.23$ for $\gamma = 1$ to $1.36$ for $\gamma = 2/3$,
and the values of $C$ ranging from $3.67$ for $\gamma = 2/7$ to $1.91$ for $\gamma = 1$. 
These plots were made for the update factor $K = 50$, which is large enough so that, for each $N$,  the errors have reached their minima.
The rest of the plots in Fig. \ref{fig:linear} are made for $N=2^{11}$.

Fig. \ref{fig:linear}(c)  indicates that the normalized maximal absolute error 
does depend on the orientation of the eigenvectors of the matrix $\Sigma$
but remains within the same order of magnitude. It changes at most by the factor of 5. 
 
Fig. \ref{fig:linear}(d) shows that the normalized maximal absolute error 
increases as $\gamma$ deviates from 1, i.e., as the condition number of the diffusion matrix increases.
The {\tt olim4vad} keeps the numerical errors reasonably small as $\gamma$ ranges from $2/7$ to $7/2$. 
This figure suggests that  {\tt olim4vad} is an appropriate method when the ratio of the maximal to the minimal eigenvalue of $A = (\sigma\sigma^{\top})^{-1}$
does not exceed 10.

Finally, Figs. \ref{fig:linear}  (e) and (f) show that the larger the ratio of the 
eigenvalues of $\Sigma$ is, the larger values of $K$ need to be used.
For each $\gamma$, the value of $\alpha$  maximizing the error was used. 
The comparison of these two figures offers a consistency test: the graphs for reciprocal values of 
$\gamma$ and the corresponding  values $\alpha$ differing
by $\pi/2$ collapse in Fig. \ref{fig:linear} as they should. 

{ The choice of the update parameter $K$ is very important for practical purposes.
This linear example shows that the error first monotonically decreases with $K$ and then stabilizes.
For our nonlinear test example in Section \ref{sec:nonlin}, the error first decreases with $K$, then stabilizes, but then starts to grow.
An optimal value for $K$ is problem-dependent and cannot be known in advance if the exact solution is unknown. 
Large values of $K$ allow to accommodate large angles between the characteristics and the gradients, 
however, they fail to account for the curvature of the characteristics which may lead to extra errors. Hence, there is a trade-off. We refer an interested reader 
to our discussion on this subject in Ref. \cite{OLIMs} (see Section 4, and, in particular, Subsection 4.3). 
From the practical point of view, it is valuable that there is a quite broad 
range of values of $K$ for which the error is nearly minimal. Based on numerical tests, we proposed the following Rule-of-Thumb \cite{OLIMs}
for choosing a reasonable value of $K$:
}
\begin{equation}
\label{Rule}
K = 10 + 4(p-7)~~\text{with}~~N = 2^p.
\end{equation}
In the present work focused on the anisotropic diffusion, we would like to check whether this rule remains applicable for this case.
Eq. \eqref{Rule} would give $K=26$ for $N=2^{11}$. Figs.  \ref{fig:linear}  (e) and (f) 
show that this estimate would be good for $2/5\le\gamma\le 5/2$, 
while it is better to increase $K$ by a factor about 1.5   for $\gamma$ closer to $0$ or to $\infty$. 
Or a larger numerical error should be tolerated.
Roughly speaking, the proposed in \cite{OLIMs} Rule-of-Thumb partially based on tests on the same linear field $J\mx$ 
with the ratio of the rotational and the potential components equal to 10, but $\Sigma = I$, 
is good for SDE \eqref{linear_test} with the ratio of the eigenvalues of the matrix 
$A =( \Sigma\Sigma^\top)^{-1}$ $\lambda_{max}/\lambda_{min}\le 10$.

\subsection{A nonlinear SDE with variable anisotropic diffusion}
\label{sec:nonlin}
A nonlinear SDE with variable anisotropic diffusion and an analytically available formula for the quasi-potential can be constructed
by designing the diffusion matrix to be a Jacobian matrix of a nonlinear variable change. The vector field should be chosen so that, 
after this variable change,
it easily decomposes to rotational and potential components. We pick $\sigma(\mx)$ for the variable change from Cartesian
 to polar coordinates 
and set up the following SDE:
\begin{equation}
\label{nonlintest}
\left[\begin{array}{c}dx_1\\dx_2\end{array}\right] = \left[\begin{array}{r} x_2g(r,\phi) +x_1f(r,\phi) \\ -x_1g(r,\phi)+x_2f(r,\phi)\end{array}\right]dt +
\sqrt{\epsilon}\left[\begin{array}{rr}r^{-1}x_1 & -x_2\\r^{-1}x_2 & x_1\end{array}\right]\left[ \begin{array}{c} dw_1\\dw_2\end{array}\right],
\end{equation}
where $r = \sqrt{x_1^2+x_2^2}$, $\phi$ is the polar angle, and the functions $f(r,\phi)$ and $g(r,\phi)$ will be specified a bit later.
{ Using the Ito formula to perform the variable change from $x_1$ and $x_2$ to $r$ and $\phi$ and neglecting the $O(\epsilon)$ 
and higher order terms, we obtain the differentials $dr$ and $d\phi$: 
$$
\left[\begin{array}{c} dr\\ d\phi\end{array}\right]
 = \left[\begin{array}{cc}r^{-1}x_1&r^{-1}x_2 \\ -r^{-2}x_2& r^{-2}x_1\end{array}\right]\left[\begin{array}{c}dx_1\\dx_2\end{array}\right].
$$}

After doing some simple algebra we obtain
\begin{equation}
\label{nonlintest1}
\left[\begin{array}{c} dr\\ d\phi\end{array}\right] = \left[\begin{array}{c} rf(r,\phi)\\ -g(r,\phi)\end{array}\right]dt +
 \sqrt{\epsilon}\left[ \begin{array}{c} dw_1\\dw_2\end{array}\right].
\end{equation}
Now it remains to pick $f(r,\phi)$ and $g(r,\phi)$ so that the deterministic term in the right-hand side of Eq. \eqref{nonlintest1}
readily decomposes into a potential and rotational components. We pick:
$$
f(r,\phi) = 1 - \frac{r^2}{9} + \frac{\sin\phi}{r^2},\quad g(r,\phi) =  \frac{\sin\phi}{r^2} - \left( 1 - \frac{r^2}{9}\right).
$$
The point $(r=3,\phi = 0)$ is an asymptotically stable equilibrium. The exact quasi-potential with respect to this equilibrium is given by
\begin{align}
U &= r^2\left(\frac{r^2}{18} - 1\right) + \frac{9}{2} + 2\left(1 -\cos\phi \right) \notag \\
&= \|\mx\|^2\left(\frac{\|\mx\|^2}{18} - 1 \right) +  \frac{9}{2} + 2\left(1 - \frac{x_1}{\|\mx\|} \right). \label{Uexact}
\end{align}
It is well-defined in any domain of the form $\mr^2\backslash \{\mx~|~\|\mx\| < r_0\}$ where $r_0>0$.
The rotational component is 
\begin{equation}
\label{rot}
\mathbf{l}(\mx) = \left[\begin{array}{c} x_2\left(\frac{\|\mx\|}{9}-1\right) + \frac{x_1x_2}{\|\mx\|^{3/2}} \\
 -x_1\left(\frac{\|\mx\|}{9}-1\right) + \frac{x_2^2}{\|\mx\|^{3/2}} \end{array}\right].
 \end{equation}
One can check that the ratio of the magnitudes of the rotational and the potential components is one
and that the eigenvalues of the matrix $A(\mx) = (\sigma(\mx)\sigma(\mx)^\top)^{-1}$ are $1$ and $\|\mx\|^{-2}$.

\begin{figure}[htbp]
\begin{center}
\centerline{
(a)\includegraphics[width = 0.4\textwidth]{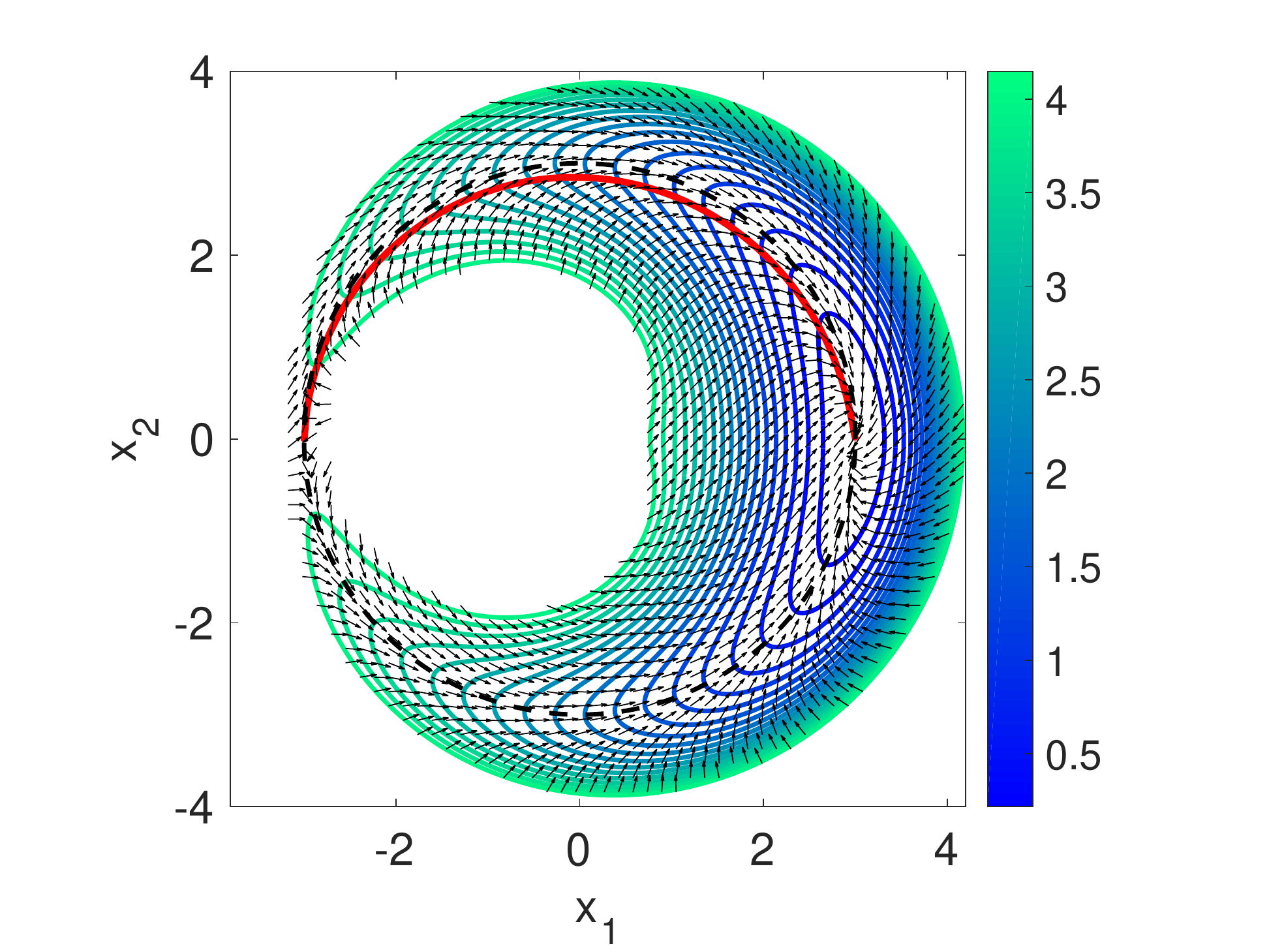}
(b)\includegraphics[width = 0.35\textwidth]{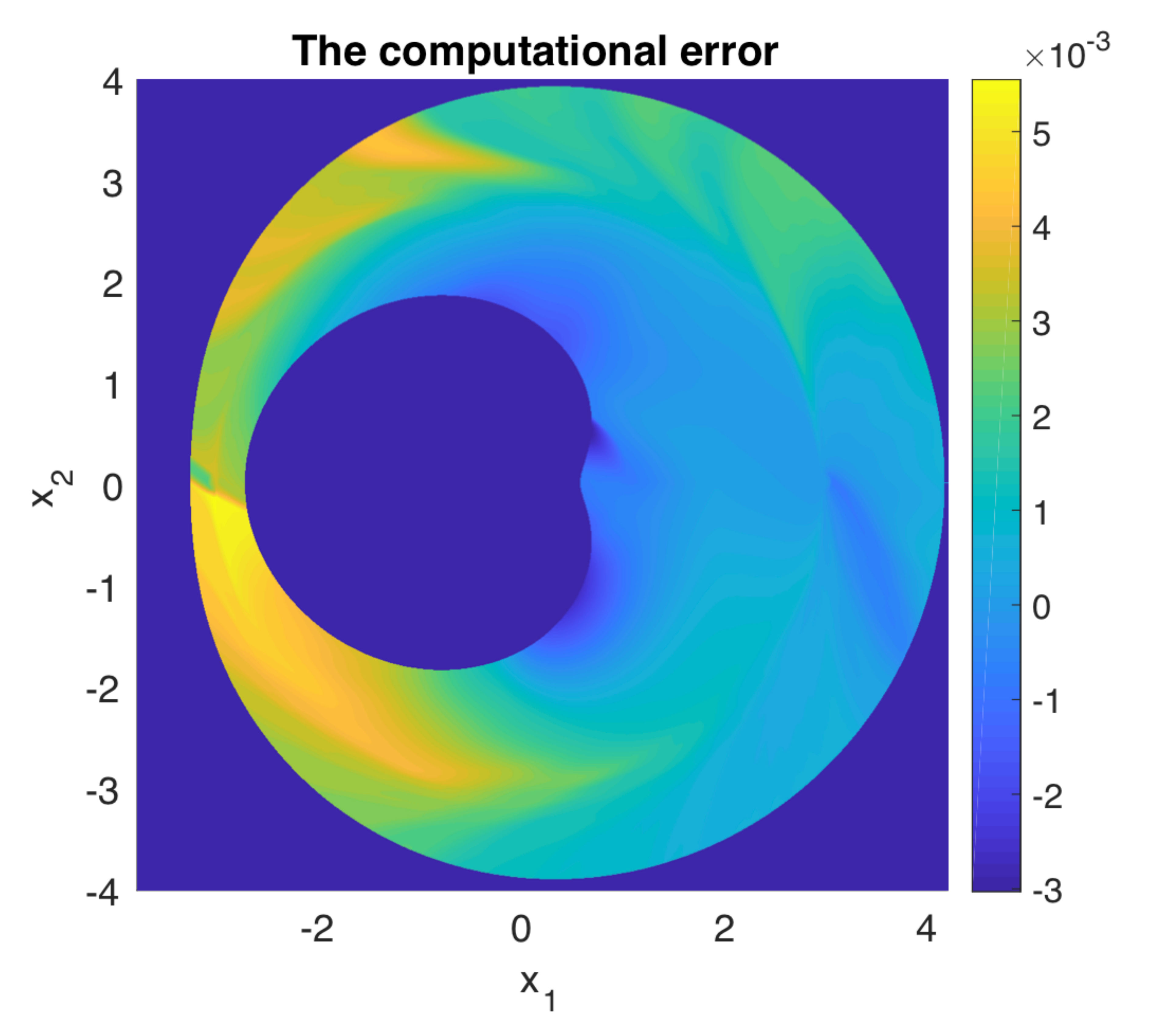}
}
\centerline{
(c)\includegraphics[width = 0.4\textwidth]{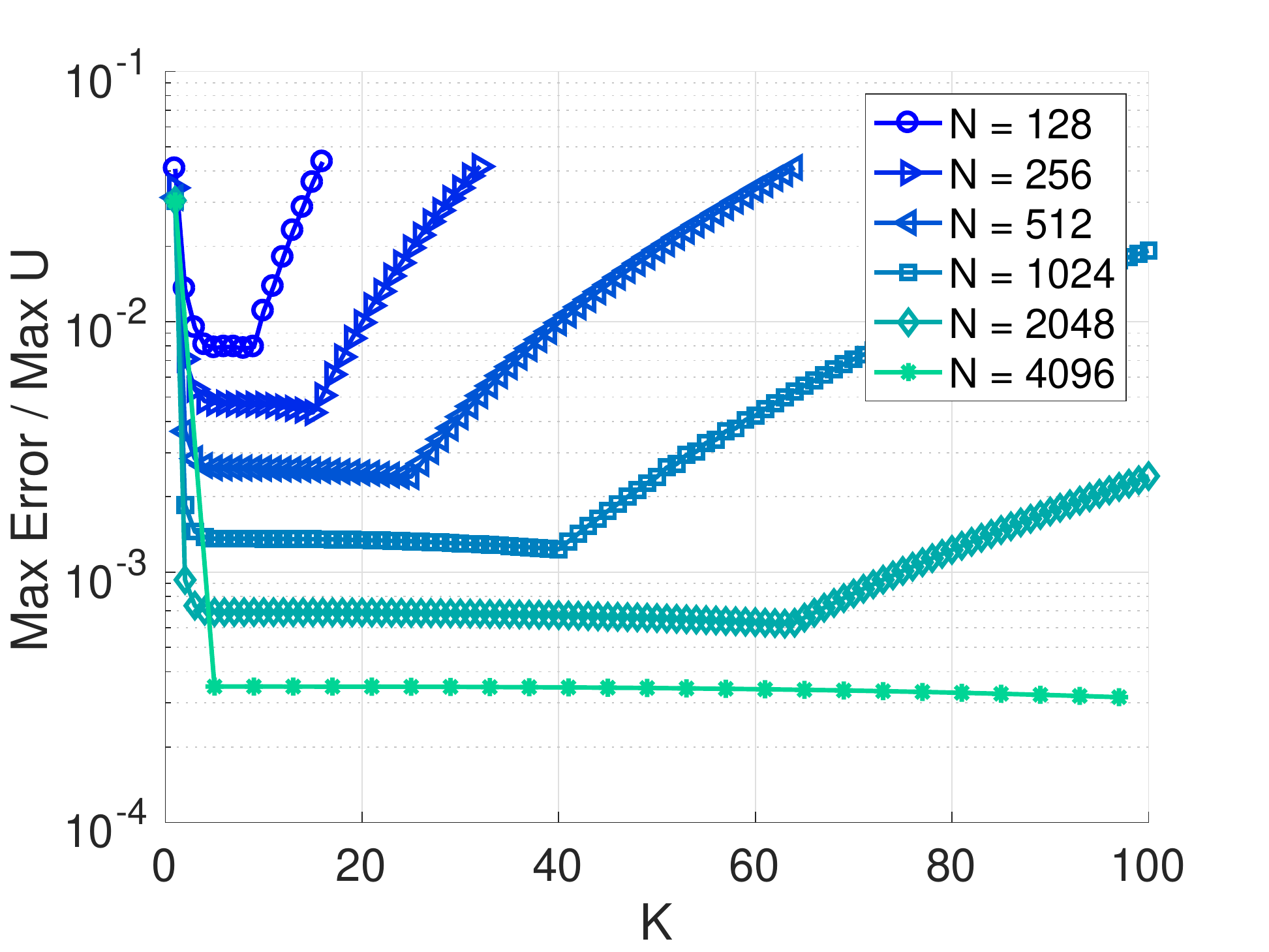}
(d)\includegraphics[width = 0.4\textwidth]{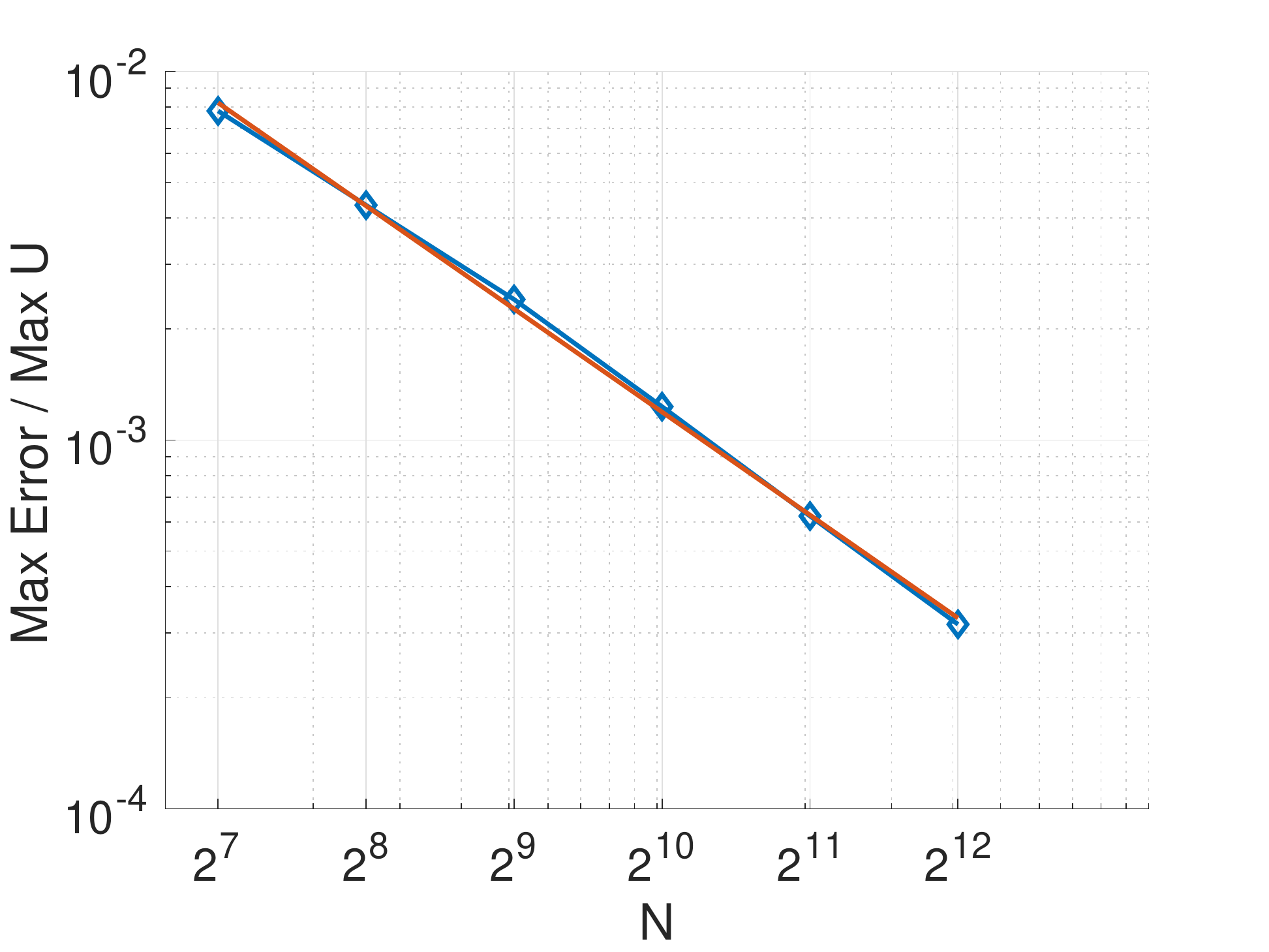}
}
\caption{
Numerical errors produced by {\tt olim4vad} for SDE \eqref{nonlintest}.
(a): The level sets of the computed quasi-potential are shown by the contour plots. 
Red curve: the MAP from the asymptotically stable equilibrium at $(3,0)$ to the saddle at $(-3,0)$.
Black arrows: the direction of the vector field. Black dashed line: the circle $r = 3$.
(b): The error plot: $U - U_{exact}$ for $N = 1024$, $K = 40$.
(c): The normalized maximal absolute error versus $K$.
(d): Blue curve with diamond markers: the normalized maximal absolute error versus $N$.
Red line: the least squares fit: $E = 0.743 N^{-0.928}$.
}
\label{fig:nonlin}
\end{center}
\end{figure}

The computational domain for this example, the square 
$$
\{-3.8\le x_1\le 4.2,~-4.0\le x_2\le 4.0\},
$$
is chosen so that the computation reaches the saddle point at $(x_1 = -3,x_2 = 0)$ { before it reaches the boundary which terminates it.
Furthermore, the origin where the quasi-potential is not defined is not reached.} 
The maximal computed value of the quasi-potential  in this domain is approximately $4.15$.
The numerical errors in Figs. \ref{fig:nonlin} (c) and (d) are normalized by it. 

The direction of the vector field, the level sets of the computed quasi-potential, and the MAP from 
the asymptotically stable equilibrium to the saddle are shown in Fig. \ref{fig:nonlin}(a).
The error plot $U - U_{exact}$ is presented in Fig. \ref{fig:nonlin}(b).
The graphs of the normalized maximal absolute error versus the update factor $K$ for $N=2^p$, $p = 7,8,\ldots,12$
 in Fig. \ref{fig:nonlin}(c) show that $K$ should not be taken too large. The Rule-of-Thumb \cite{OLIMs} suggesting  $K = 10 + 4(p-7)$
 gives a good set of values of $K$ for this example. The graph of the normalized maximal absolute error vs $N$ and the least squares fit  $E = CN^{-p}$
 are shown in Fig. \ref{fig:nonlin}(d). The least squares fit gives: $C = 0.743$, $p = 0.928$.


\section{A demo: the Maier-Stein model with anisotropic diffusion}
\label{sec:MS}
In this Section, we  demonstrate the effects of anisotropy on the well-known Maier-Stein model \cite{MS1993,MS1996}.

The original Maier-Stein model is given by SDE \eqref{sde1} with $\sigma(\mx)\equiv I$ (the identity matrix) and the vector field
\begin{equation}
\label{MSfield}
\mb(\mx) = \left[\begin{array}{c} x_1 - x_1^3 - 10x_1x_2^2 \\ -(1.0 + x_1^2)x_2  \end{array}\right].
\end{equation}
It has two asymptotically stable equilibria at $(-1,0)$ and $(1,0)$ and a saddle point at the origin. The vector field is symmetric with respect to
the $x_1$-axis and the $x_2$-axis.
The quasi-potential with respect to either equilibrium is nondifferentiable along a segment of the $x_1$-axis \cite{quasi}.

\begin{figure}
    \centering
    \begin{subfigure}[b]{0.4\textwidth}
        \includegraphics[width=\textwidth]{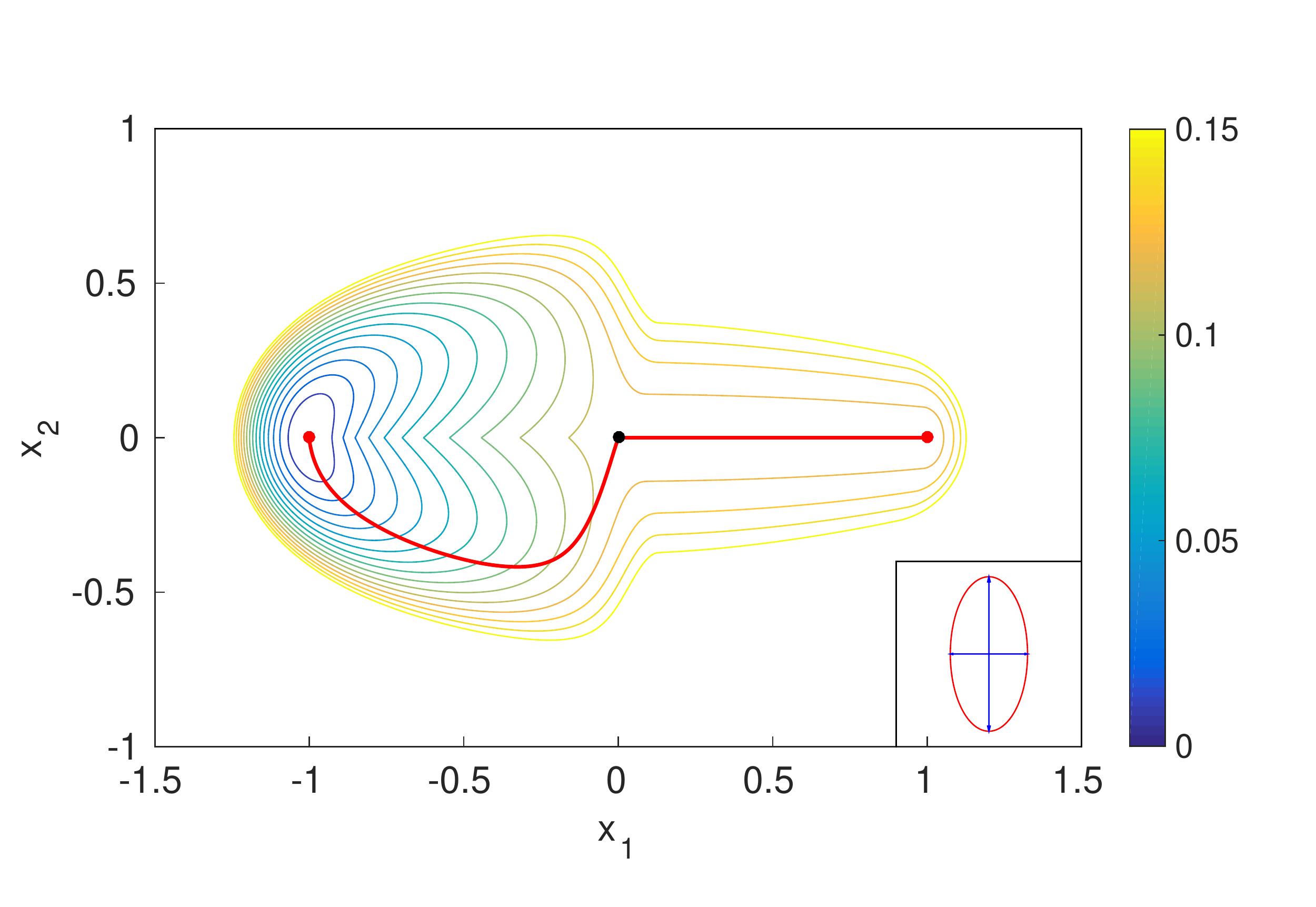}
        \caption{$\alpha = 0$}
    \end{subfigure}
    ~ 
    \begin{subfigure}[b]{0.4\textwidth}
        \includegraphics[width=\textwidth]{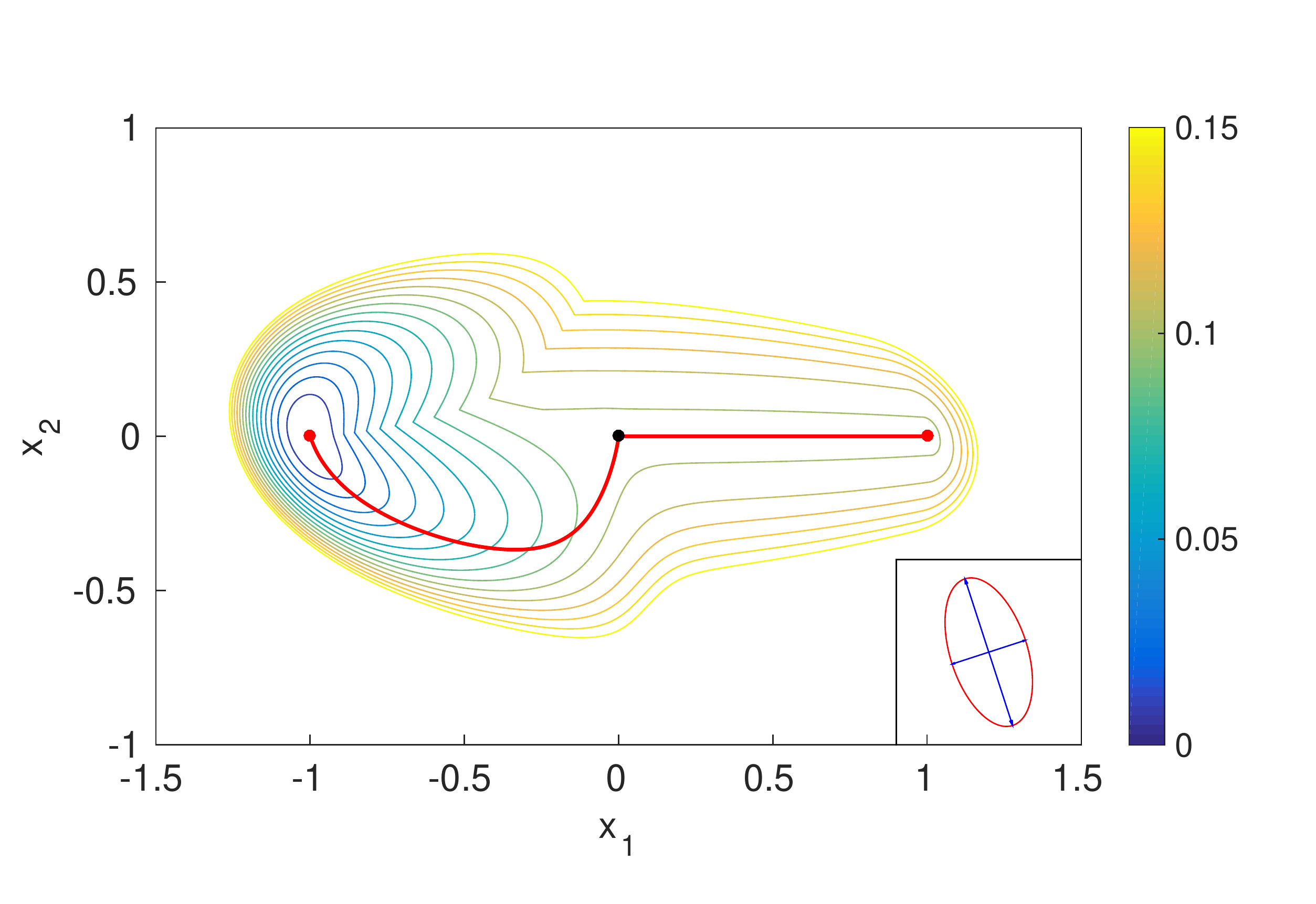}
        \caption{$\alpha = \pi/10$}
    \end{subfigure}
    ~ 
    \begin{subfigure}[b]{0.4\textwidth}
        \includegraphics[width=\textwidth]{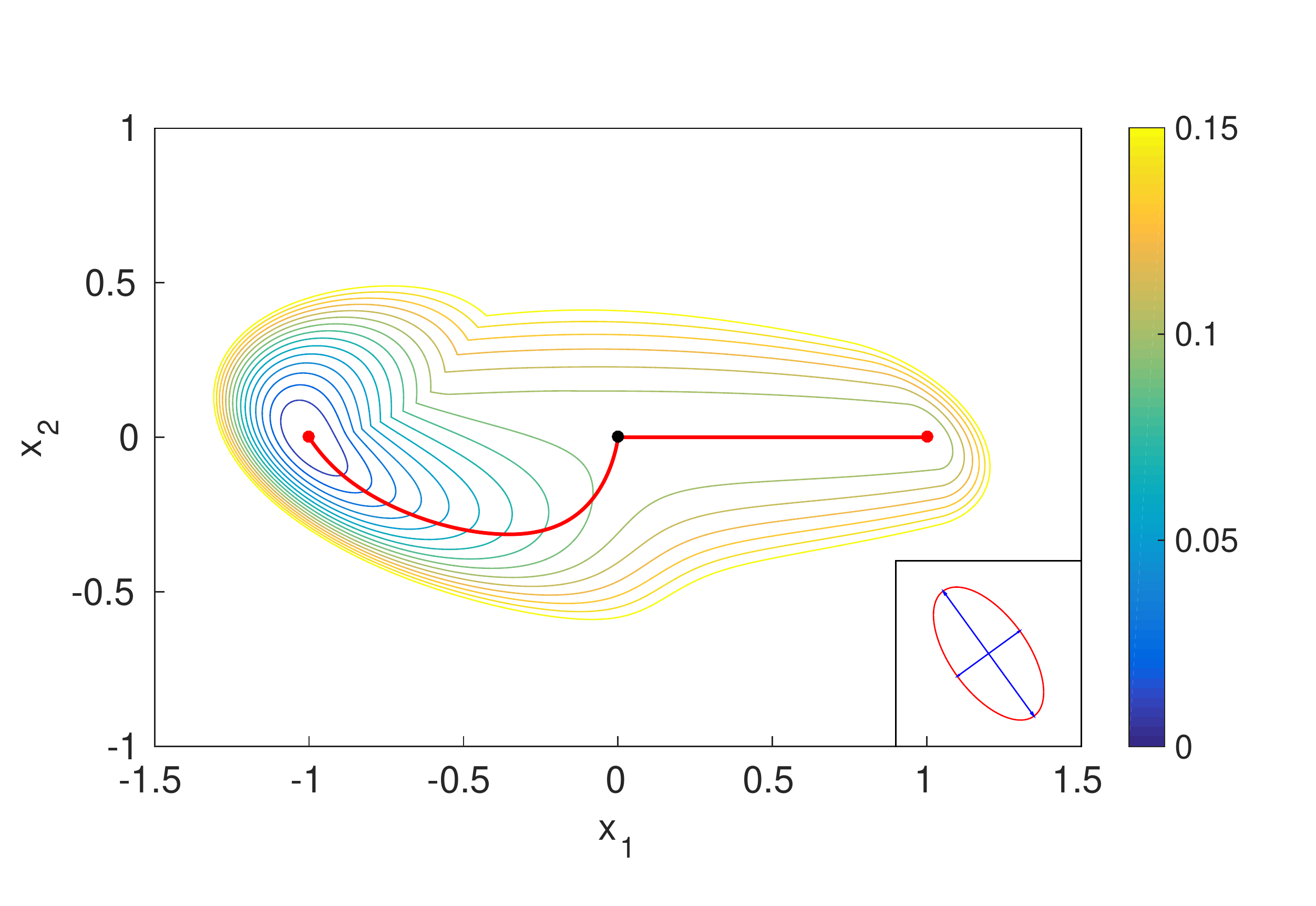}
        \caption{$\alpha = \pi/5$}
    \end{subfigure}
        \begin{subfigure}[b]{0.4\textwidth}
        \includegraphics[width=\textwidth]{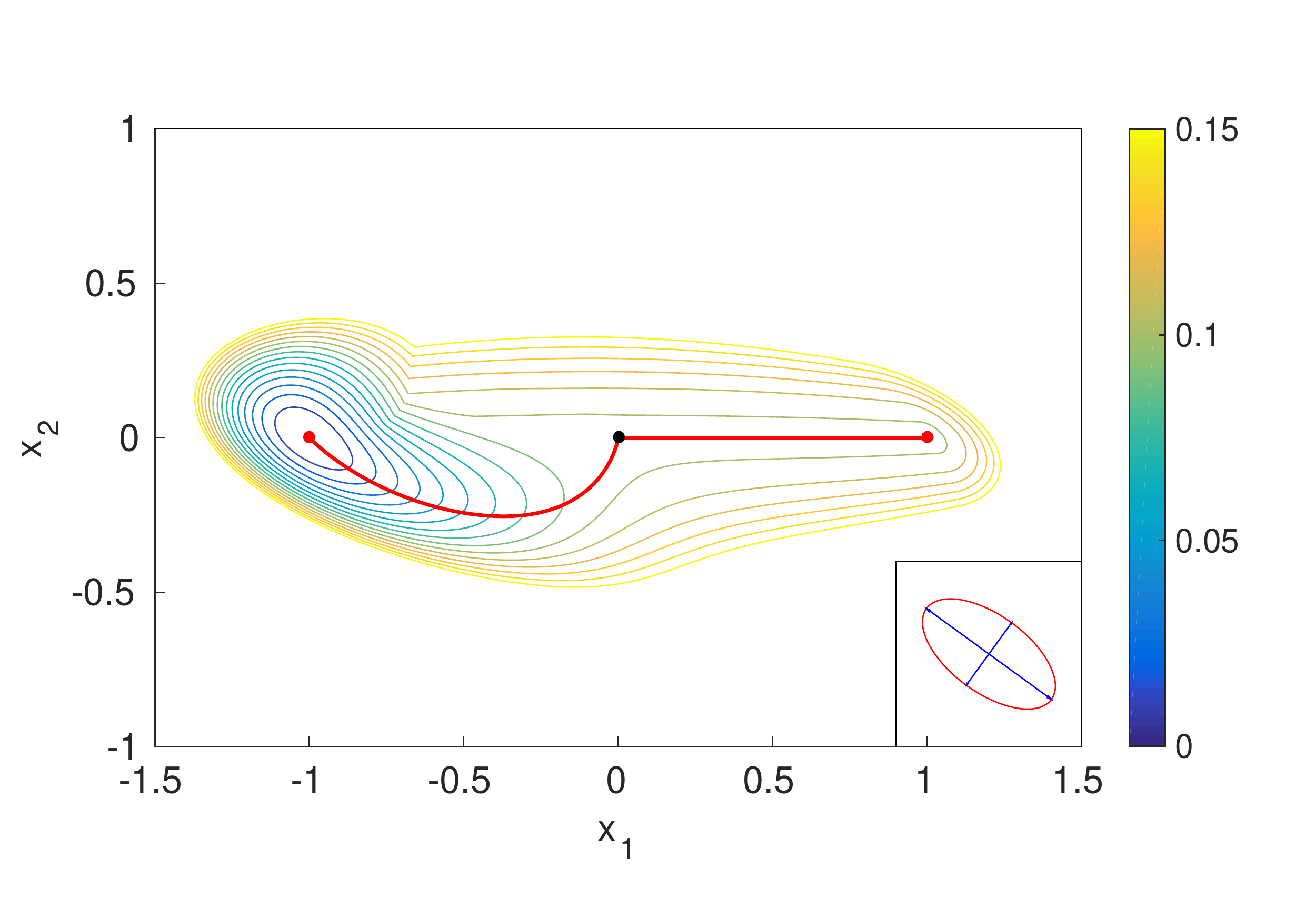}
        \caption{$\alpha = 3\pi/10$}
    \end{subfigure}
    ~ 
    \begin{subfigure}[b]{0.4\textwidth}
        \includegraphics[width=\textwidth]{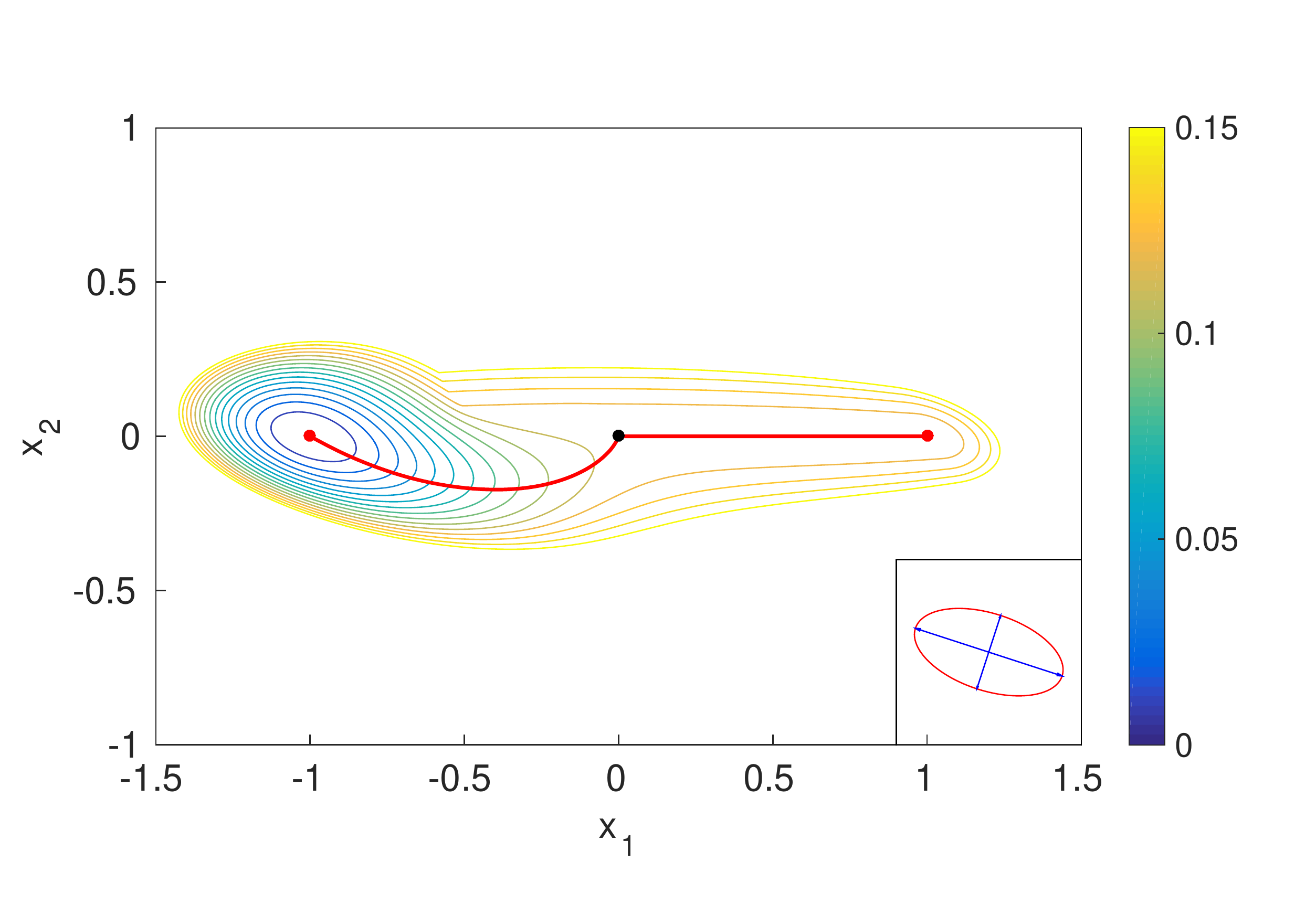}
        \caption{$\alpha = 2\pi/5$}
    \end{subfigure}
    ~ 
    \begin{subfigure}[b]{0.4\textwidth}
        \includegraphics[width=\textwidth]{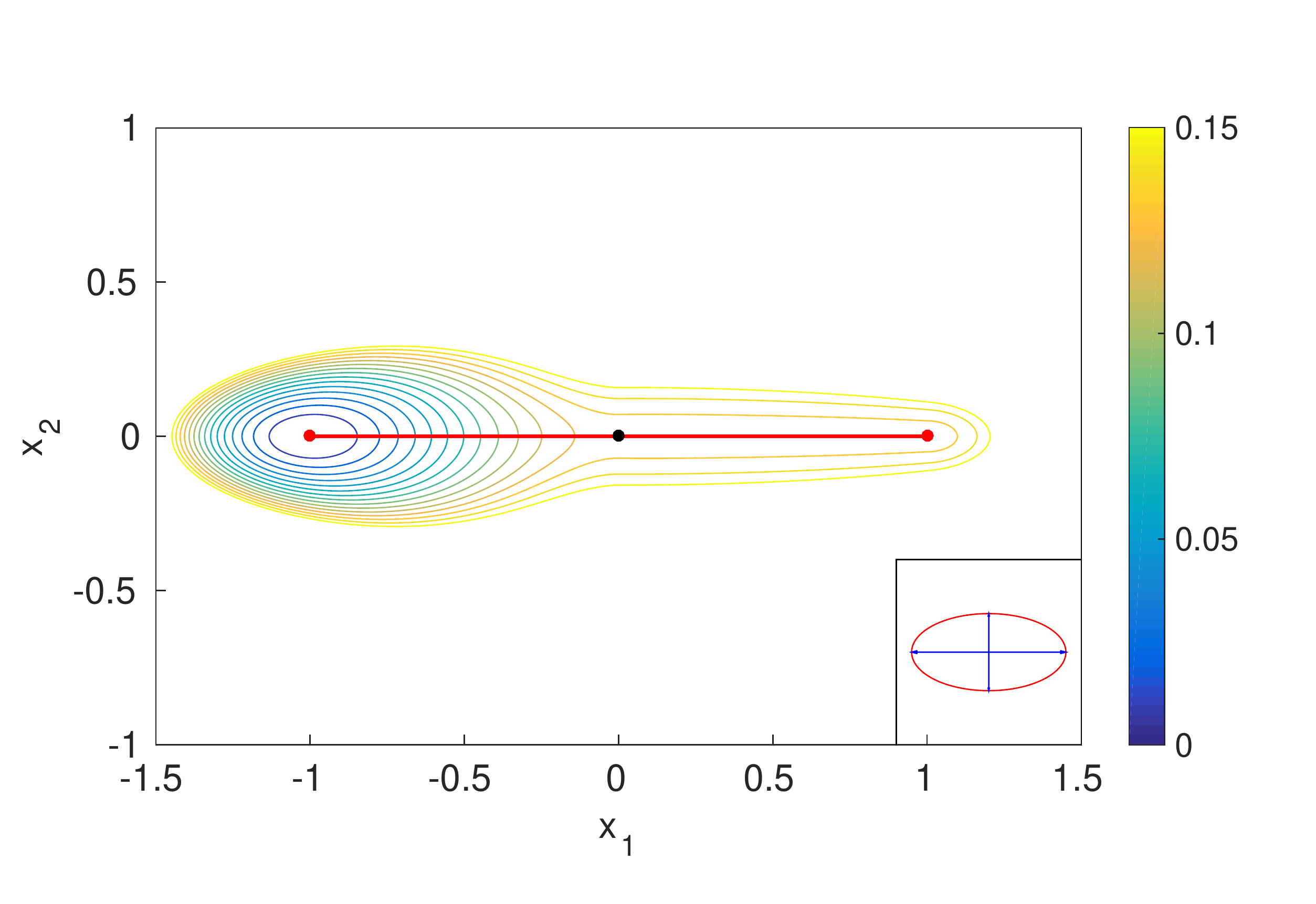}
        \caption{$\alpha = \pi/2$}
    \end{subfigure}
    \caption{The red points are the two stable equilibria of the Maier-Stein model and the 
    black point is the saddle. The figures show the level sets of the quasi-potential and 
    the MAPs (red curves) for the Maier-Stein model for the diffusion matrices 
    corresponding to $\gamma = 2$ and $\alpha \in \lbrace 0, \pi/10, \pi/5, 3\pi/10,2\pi/5,\pi/2\rbrace$ for $N=2048$ and $K=25.$ The lower right corner in each figure shows the level of anisotropy.}
\label{fig:MS}
\end{figure}

In this work, we have applied {\tt olim4vad} to the Maier-Stein model with a family of constant diffusion matrices of the form
\begin{equation}
\label{MSsigma}
\Sigma =\left[\begin{array}{cc} \cos \alpha & -\sin \alpha \\ \sin \alpha & \cos \alpha \end{array}\right] \left[\begin{array}{cc} 1 & 0
 \\ 0 & \gamma \end{array}\right] \left[\begin{array}{cc} \cos \alpha & -\sin \alpha \\ \sin \alpha & \cos \alpha \end{array}\right]^{-1}.
\end{equation}
We have set $\gamma =2$ and run {\tt olim4vad} for the following set of values of the parameter $\alpha$
determining the orientation of the eigenvectors of $\Sigma$:
\begin{equation}
\label{MS_alpha}
\alpha \in \left\{ 0, \frac{\pi}{10}, \frac{\pi}{5}, \frac{3\pi}{10}, \frac{2\pi}{5},\frac{\pi}{2}\right\}.
\end{equation}
It suffices to consider just this set of values of $\alpha$ due to the symmetry with respect to the $x_1$-axis.
The computations were conducted in the domain $[-2,2]^2$ starting at the equilibrium $(-1,0).$ The update factor $K = 25$ and $N=2048$ were chosen. 
The level sets of the computed quasi-potential and the MAPs are shown in Fig. \ref{fig:MS}.
The lower right corner in each figure displays the 
orientation of the anisotropy. We observe that the quasi-potential and the 
MAPs change significantly as the orientation of anisotropy changes. 
{ The symmetry of the original Maier-Stein model is broken unless the eigenvectors of $\Sigma$ are parallel to the coordinate axes.

The computation of the quasi-potential in 2D is cheap. For example,  the computation for the Maier-Stein model on 
the $2048\times2048$ mesh with $K = 25$ takes about 16 seconds.
In comparison with path-based methods, the computation of the quasi-potential gives more information about the 
asymptotic behavior of the system. The MAP is readily obtained. The equilibrium probability density is  available up to the prefactor:
$$
\mu(\mx)\asymp e^{-U(\mx)/\epsilon},
$$
where $\mx$ belongs to the sublevel set of the quasi-potential passing through the saddle.
The visualization of the quasi-potential allows one to understand the asymptotic behavior of the system just at a glance.  
}

The quasi-potential of the Maier-Stein is non-differentiable at the saddle point $(0,0)$. Therefore, the Bouchet-Reygner 
formula \eqref{TransTime}  for the sharp estimate for the expected exit time from the basin of attraction is not applicable for this example.

\section{An application to a genetic toggle model for Lambda Phage.}
\label{sec:GT}
In this Section, we apply {\tt olim4vad} to a genetic toggle model for Lambda Phage  
\cite{GTbook,Lei,SheaAckers,AurellBrown,Aurell,Wang} 
and use the Bouchet-Reygner formula \eqref{TransTime} to
calculate the expected transition time from the lysogenic to the lytic state.

\begin{figure}
\centering
\includegraphics[width = 0.65\textwidth]{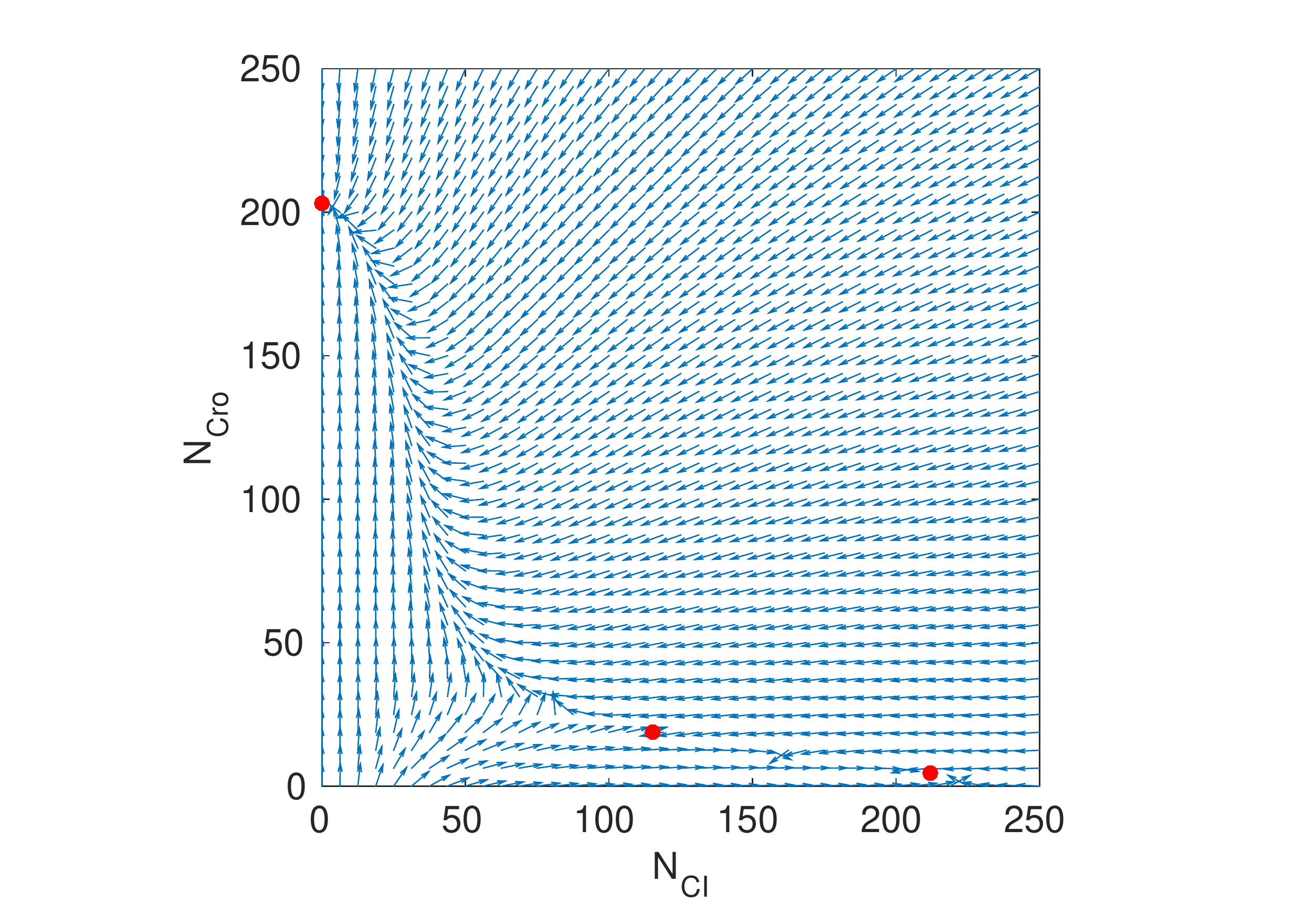}
\caption{The vector field for the genetic toggle problem. 
Red dots mark two stable equilibria $(0.1654,203.0115)$ (the lytic state) 
and $(212,4.5)$ (the lysogenic state), and a transition state at $(115.0625,18.6875)$.
} 
\label{fig:VectorFieldGT}
\end{figure}

Lambda Phage is a virus which can infect the bacterium Escherichia coli. 
Once a bacterial cell is infected, the virus can adopt two different lifecycles:  
lysogenic (in which the virus 
remains inside the bacterial cell but does not replicate) or lytic 
(in which the bacterial cell dies as the virus replicates). 
Which one of these two cases will happen 
depends on the gene expression \cite{GTbook}. This problem has been studied by many researchers. 
Shea and Ackers 
\cite{SheaAckers} developed a quantitative physical-chemical model for this problem.  Aurell et. al. \cite{AurellBrown} 
developed a stochastic model to determine the relative stability of the two states based on the molecular interactions. 
The choice between lysis and lysogeny is mainly governed by two regulatory proteins, CI and Cro, 
which are encoded by viral genes {\it ci} and {\it cro} respectively. The sequence between these two 
genes comprises two promoters (${\rm P_{RM}}$ and ${\rm P_R}$) and an operator ${\rm O_R}$ 
(see figure 1.4 in \cite{GTbook}). These two regulatory proteins, along with RNA polymerase, promoters 
${\rm  P_{RM}}$ and ${\rm P_R}$, and an operator ${\rm O_R}$, constitute the genetic switch which 
decides between the two possible pathways a viral DNA may undertake in the host cell. 
The sources of noise for this problem are discussed by Lei. et. al. \cite{Lei}. 
Aurell et al \cite{Aurell} proposed the following SDE to model the counts $N_{CI}$ and $N_{Cro}$ 
of the key proteins CI and Cro:
\begin{align}
dN_{CI} =& \left[ S_{CI} f_{CI}(N_{CI},N_{Cro})-N_{CI}/\tau_{CI}\right]dt+g_{CI} dW^{CI},\notag\\
dN_{Cro} =& \left[S_{Cro} f_{Cro}(N_{CI},N_{Cro})-N_{Cro}/\tau_{Cro}\right]dt+g_{Cro} dW^{Cro}\label{GTSDE}
\end{align}
The diffusion matrix $\sigma(N_{CI},N_{Cro}) = {\rm diag}(g_{CI},g_{Cro})$ in Eq. \eqref{GTSDE} is diagonal and position-dependent: 
\begin{align}
g_{CI} = \sqrt{S^2_{CI}f_{CI}+N_{CI}/\tau_{CI}},~~g_{Cro} = \sqrt{S^2_{Cro}f_{Cro}+N_{Cro}/\tau_{Cro}}.\label{gCICro}
\end{align}
 The functions $f_{CI},$ $f_{Cro}$ are quite complicated and placed in Appendix C.
All parameters for Eq. \eqref{GTSDE} are also found in Appendix C.
The direction of the vector field in Eq. \eqref{GTSDE}
is plotted  in Fig. \ref{fig:VectorFieldGT}. 
The deterministic system corresponding to the SDE \eqref{GTSDE} is bistable with two asymptotically 
stable equilibria at $(0.1654,203.0115)$ (the lytic state) and $(212,4.5)$ (the lysogenic state).
By computing the quasi-potential with respect to the lysogenic state, we found the transition state located at 
$(115.0625,18.6875)$ (see Fig. \ref{fig:VectorFieldGT}): after reaching this point, 
the quasi-potential remains constant along the trajectory going to the lytic state.
This location is different from the one suggested in \cite{Aurell}.
Note that the maximum and minimum magnitude of the vector field are $\approx 6.0$ 
and $\approx 1.1\times10^{-4}$ respectively. 
The vector field has very small magnitude in the transition channel between the lysogenic and the lytic states.

\begin{figure}
\centerline{
\includegraphics[width=0.8\textwidth]{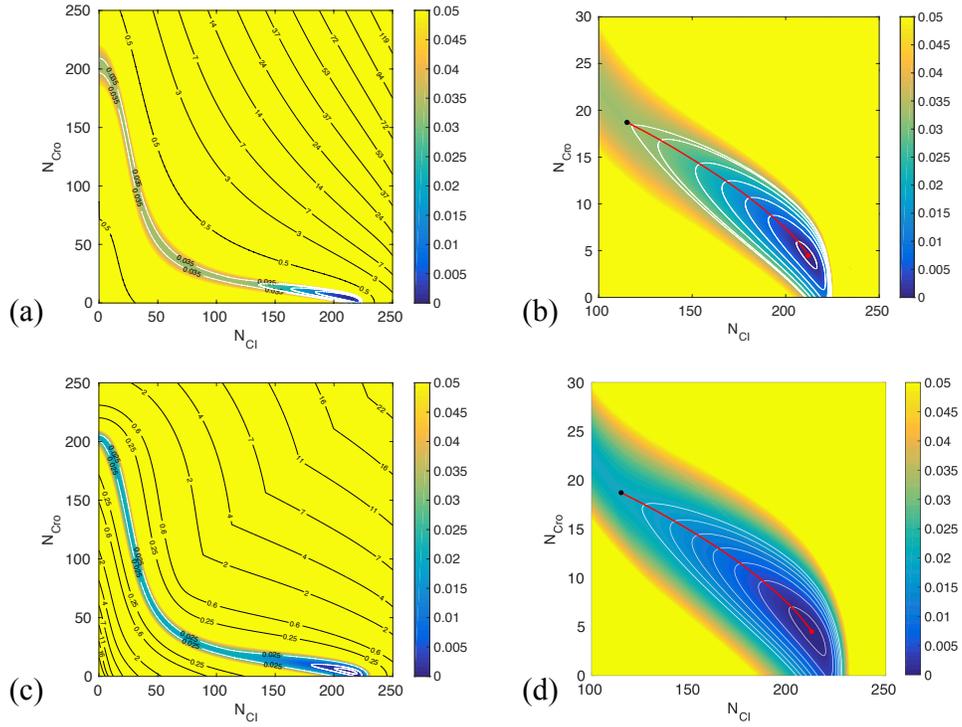}
}
\caption{
 (a):  The quasi-potential with respect to the lysogenic state and the MAP  connecting the lysogenic state and the transition state
 computed for the SDE \eqref{GTSDE}.
 (b): A zoom-in of the MAP and the region around it.
 (c) and (d): The same for the diffusion matrix in SDE \eqref{GTSDE}
 replaced with the identity matrix.
}
\label{fig:GT}
\end{figure}

The computation of the quasi-potential has been done in the region $[0,250]\times[0,250]$  on $2048\times 2048$ mesh with $K=25$.
Since the rotational component in this system is small, 
the computation has continued throughout the whole square $[0,250]^2$.
The resulting quasi-potential and the MAP from the lysogenic state to the transition state are shown 
in Fig. \ref{fig:GT} (a) and (b).

For comparison, we have also computed the quasi-potential with respect to the lysogenic state for uniform isotropic diffusion, i.e., $\sigma(N_{CI},N_{Cro})\equiv I$.
We observe that the two MAPs differ insignifcantly, however the quasi-potentials differ notably. 
Our computation of the MAPs locates the saddle point which is otherwise very difficult to find as the vector field is close to zero in a whole region
and the formula for the vector field is very complicated.
The transition rates computed for both choices of the diffusion matrix using Bouchet-Reygner formula \eqref{TransTime}
are shown in
Table \ref{Table:TransRate}.  { The computation of the quasi-potential has also shown that the transition channel from the lysogenic to the lytic state 
is very narrow. This means that, even if the noise term is not small, the transition process is effectively bound to this channel. }

\begin{table}
\caption{The transition rate of the genetic toggle system from the lysogenic state to the lytic state.}
\label{Table:TransRate}     
\centering
\begin{tabular}{lllll}
\hline\noalign{\smallskip}
Type of Diffusion Matrix & Transition Rate $s^{-1}$ \\
\noalign{\smallskip}\hline\noalign{\smallskip}
Identity  &  $3.76072\times 10^{-6}$\\
Diagonal &  $4.29250\times 10^{-6}$\\
%
\noalign{\smallskip}\hline
\end{tabular}
\end{table}


\section{Conclusion}
\label{sec:conclusion}
In this work,  the Ordered Line Integral Method (OLIM) with the midpoint quadrature rule is extended
for computing the quasi-potential for stochastic differential equations with position-dependent and anisotropic diffusion.
We named the resulting method {\tt olim4vad}. 
Our C code {\tt OLIM4VAD.c} is available at M. Cameron's website \cite{Mweb}.

{\tt olim4vad} has been tested on two examples where the quasi-potential can be found analytically:
a linear SDE with a family of constant diffusion matrices, and a nonlinear SDE with a position-dependent diffusion matrix.
Our analysis of the dependence of the numerical 
error on $N$ (the mesh is $N\times N$), the update factor $K$, 
the ratio of the eigenvalues, and the direction of the eigenvectors of the diffusion matrix 
suggests that the Rule-of-Thumb for choosing $K$
for a given $N$ proposed in \cite{OLIMs} 
is good for the case where the ratio $\lambda_{\max}/\lambda_{\min}$ of the matrix $A=(\sigma\sigma^\top)^{-1}$ does not exceed 10.

We have demonstrated the effects of anisotropy 
on the Maier-Stein model with a family of constant diffusion matrices 
with the ratio of eigenvalues equal to 2 and various directions of the 
eigenvectors. 
The quasi-potential level sets and the MAPs 
get significantly deformed as the orientation of the anisotropy changes.

Finally, we have applied {\tt olim4vad} to the genetic toggle model for Lambda Phage \cite{Aurell},
an SDE with a position-dependent diffusion matrix. We have found the transition state between the lysogenic and lytic states.
We have also found the escape  rate from the lysogenic state using the Bouchet-Reygner formula \cite{Bouchet}.

\section{Acknowledgements}
This work was partially supported by M. Cameron's NSF grant DMS1554907.


 \appendix
\setcounter{equation}{0}
\renewcommand{\theequation}{A-\arabic{equation}}

{
 \section*{Appendix A. A derivation of the geometric action}
 \label{sec:AppA}
 
The minimization with respect to the travel time $T$ for the Freidlin-Wentzell action \eqref{FWA}
can be performed analytically as follows \cite{hey2,FW}:
\begin{align}
S_T(\phi) &= \frac{1}{2}\int_0^T\|\dot{\phi} - \mb(\phi)\|^2_{A(\phi)}dt = \frac{1}{2}\int_0^T(\dot{\phi}- \mb(\phi))^\top A(\phi) (\dot{\phi}- \mb(\phi))dt \notag \\
&=\frac{1}{2}\int_0^T \left(\|\dot{\phi}\|^2_{A(\phi)} - 2\langle\dot{\phi},\mb(\phi)\rangle_{A(\phi)} +  \|\mb(\phi)\|^2_{A(\phi)}\right)dt \notag \\
&\ge  \frac{1}{2}\int_0^T \left(2\|\dot{\phi}\|_{A(\phi)} \|\mb(\phi)\|_{A(\phi)} - 2\langle\dot{\phi},\mb(\phi)\rangle_{A(\phi)}\right)dt \label{derg} \\
& = \int_0^T\left(\|\dot{\phi}\|_{A(\phi)} \|\mb(\phi)\|_{A(\phi)} - \langle\dot{\phi},\mb(\phi)\rangle_{A(\phi)}\right)dt.\notag
\end{align}
The inequality \eqref{derg} becomes an equality if and only if $\|\mb(\phi)\|_{A(\phi)}=\|\dot{\phi}\|_{A(\phi)}.$ 
We take a reparametrization $\psi$ of the path $\phi$ such that $\|\mb(\psi)\|_{A(\psi)}=\|\dot{\psi}\|_{A(\psi)}.$
Then
\begin{equation}
\label{derg2}
S_T(\phi)\ge S(\psi) =  \int_0^{T_{\psi}}\left(\|\dot{\psi}\|_{A(\psi)} \|\mb(\psi)\|_{A(\psi)} - \langle\dot{\psi},\mb(\psi)\rangle_{A(\psi)}\right)dt.
\end{equation}
$S(\psi)$ is known as the geometric action. $T_{\psi}$ is infinite if at least one endpoint of the path is an equilibrium. 
The integral in Eq. \eqref{derg2}  is independent of the parametrization of the path $\psi$. 
We choose the arc length parametrization to obtain 
\begin{equation}
\label{GeoAction}
S(\psi) =  \int_0^{L}\left(\|\mb(\psi)\|_{A(\psi)} \|\psi_s\|_{A(\psi)} - \langle\psi_s,\mb(\psi)\rangle_{A(\psi)}\right)ds, 
\end{equation}
where $\psi_s$ is the derivative of the path $\psi$ with respect to its arc length and $L$ is the length of the path $\psi.$

}

%
%

 \appendix
\setcounter{equation}{0}
\renewcommand{\theequation}{B-\arabic{equation}}

 \section*{Appendix B. A derivation of the Hamilton-Jacobi equation for  the quasi-potential}
 \label{sec:AppB}

Let us set $\epsilon >0$ and pick such a parametrization of the path $\psi$ that  $\|\dot{\psi}\|_{A(\psi)} = 1$. 
Using Bellman's optimality principle \cite{bellman} and  Taylor expansion of $U$, we obtain
\begin{align*}
U(\mx)=\inf_{\|\dot{\psi}\|_{A(\psi)}=1}\left\{ \int_{0}^{\epsilon} \left( \|\mb (\psi)\|_{A(\psi)} -\mb(\psi)^{\top} A(\psi) \dot{\psi} \right) ds +
U\left( \mx-\int_{0}^{\epsilon} \dot{\psi} ds\right)      \right\}\\
= \inf_{\|\dot{\psi}\|_{A(\psi)}=1}\left\{ \epsilon \left( \|\mb (\psi)\|_{A(\psi)} -\mb(\psi)^{\top} A(\psi) \dot{\psi} - \nabla U(\mx)^{\top} \dot{\psi} \right)  +
U(\mx) +o(\epsilon^2) \right\}.
\end{align*}
Cancelling $U(\mx)$ on both sides and dividing by $\epsilon$ we get
\begin{align*}
0=\inf_{\|\dot{\psi}\|_{A(\psi)}=1}\left\{  \|\mb (\psi)\|_{A(\psi)} -\mb(\psi)^{\top} A(\psi) \dot{\psi} -\nabla U(\mx)^{\top} \dot{\psi} +  o(\epsilon)    \right\}.
\end{align*}
Now taking the limit as $\epsilon \rightarrow 0,$ we obtain
\begin{align}
\label{HJder}
\inf_{\|\dot{\psi}\|_{A(\mx)}=1}\left\{  \|\mb (\mx)\|_{A(\mx)} - \left(\mb(\mx)+  A(\mx)^{-1} \nabla U(\mx) \right)^{\top} A(\mx) \dot{\psi}  \right\}=0.
\end{align}
The infimum is attained when the term
\begin{equation}
\label{max} 
\left(\mb(\mx)+  A(\mx)^{-1} \nabla U(\mx) \right)^{\top} A(\mx) \dot{\psi} = \langle\mb(\mx)+A(\mx)^{-1}\nabla U(\mx),\dot{\psi}\rangle_{A}
\end{equation}
is maximal. Eq. \eqref{max} implies that  the maximizing path $\psi$ is such that
\begin{equation}
\label{psimax}
 \dot{\psi}=\displaystyle\frac{\mb(\mx) + {A(\mx)}^{-1} \nabla U(\mx)}{\| \mb(\mx) + {A(\mx)}^{-1} \nabla U(\mx) \|_{A(\mx)}}.
\end{equation}
Substituting $\dot{\psi}$ in equation \eqref{HJder}, we get
\begin{equation}
\label{HJder1}
\|\mb(\mx)\|_{A(\mx)}  = \| \mb(\mx) + {A(\mx)}^{-1} \nabla U(\mx) \|_{A(\mx)}.
\end{equation}
Taking squares of both sides of Eq. \eqref{HJder} and canceling $\|\mb(\mx)\|^2_{A(\mx)}$ we obtain the 
desired Hamilton-Jacobi equation:
\begin{equation}
\nabla U(\mx)^{\top} {A(\mx)}^{-1} \nabla U(\mx)+ 2 \mb(\mx)^{\top} \nabla U(\mx)=0.
\end{equation}

%
%

\appendix
\setcounter{equation}{0}
\renewcommand{\theequation}{C-\arabic{equation}}
\section*{Appendix C. The vector field for genetic toggle switch}
In this Appendix,  we provide all details gathered from \cite{Aurell,Wang} that are necessary for programming the genetic toggle model for Lambda Phage  given by SDE \eqref{GTSDE}.
The functions $f_{CI}$ and $f_{Cro}$  are 
\begin{align}
f_{CI}  &= R_{RM}(P_{010}+P_{011}+P_{012})+R^u_{RM}(P_{000}+P_{001}+P_{002}+P_{020}+P_{021}+P_{022}), \label{fCI}\\
f_{Cro} & = R_R(P_{020}+P_{021}+P_{022}). \label{fCro}
\end{align}
The quantities $P_s$ denote the probabilities of states $s$. 
State $s$ indicates whether the three operator sites on the DNA are free, or occupied by $CI$, 
or occupied by $Cro$ (encoded by 0, or 1, or 2 respectively). 
The total number of states is $27$.
The probability $P_s$ of state $s$ is given by
\begin{equation}
\label{ps}
P_s = \frac{[CI]^{i_s}[Cro]^{j_s}e^{-G(s)/{RT}}}{ \sum_{s} [CI]^{i_s}[Cro]^{j_s}e^{-G(s)/{RT}}},
\end{equation}
where $[CI]$ and $[Cro]$ are the concentrations of $CI$ and $Cro$ dimers respectively, $G(s)$ is the binding energy for state $s$ listed in Table \ref{Table:BindEnergy}, 
$i_s$ and $j_s$ are the numbers of $CI$ and $Cro$ molecules bound to operator sites in state $s$.
The concentrations of $CI$ and $Cro$ dimers relate to the concentrations $[N_{CI}]$ and $[N_{Cro}]$ 
of the corresponding monomers via \cite{Wang}\footnotemark[1]
\footnotetext[1]{There is an error in the expressions of $[CI]$ and $[Cro]$ in \cite{Wang} in the sign in front of the square roots. 
}
\begin{align}
[CI] &= \frac{1}{2}[N_{CI}] + \frac{1}{8}e^{\Delta G_{CI}/RT}- \left([N_{CI}] \frac{1}{8}e^{\Delta G_{CI}/RT}+\frac{1}{64}e^{2\Delta G_{CI}/RT}\right)^{1/2}, \label{CIdimer}\\
[Cro] &= \frac{1}{2}[N_{Cro}] + \frac{1}{8}e^{\Delta G_{Cro}/RT}- \left([N_{Cro}] \frac{1}{8}e^{\Delta G_{Cro}/RT}+\frac{1}{64}e^{2\Delta G_{Cro}/RT}\right)^{1/2}.
\label{Crodimer}
\end{align}
 In order to convert the numbers of $CI$ and $Cro$ into concentrations of their monomers $[N_{CI}]$ and $[N_{Cro}]$, 
one needs the effective bacterial volume, which is taken as $2\times 10^{-15}$ \cite{Aurell} 
\begin{equation}
\label{conc}
[N_{CI}] = \frac{N_{CI}}{2\times 10^{-15} N_A},~~~[N_{Cro}] = \frac{N_{Cro}}{2\times 10^{-15} N_A},
\end{equation}
where $N_A = 6.022140857 \times 10^{23}$ is the Avogadro number.
The rest of the necessary parameters for SDE \eqref{GTSDE} are listed in Table \ref{Table:param}.

\begin{table}
\caption{For each of the $27$ states, the binding energy $G(s),$ for state $s$ and the numbers of $CI$ and $Cro$ molecules, that is, $i_s$ and $j_s,$ bound to operator sites in state $s$ \cite{Aurell}.}
\label{Table:BindEnergy}     
\centering
\begin{tabular}{lllll}
\hline\noalign{\smallskip}
$s$ & $i_s$ & $j_s$ & $G(s)$ \\
\noalign{\smallskip}\hline\noalign{\smallskip}
$000$ & $0$ & $0$ & $~0$ \\
$001$ & $1$ & $0$ & $-12.5$ \\
$010$ & $1$ & $0$ & $-10.5$ \\
$100$ & $1$ & $0$ & $-9.5$ \\
$011$ & $2$ & $0$ & $-25.7$ \\
$101$ & $2$ & $0$ & $-22.0$ \\
$110$ & $2$ & $0$ & $-22.9$ \\
$111$ & $3$ & $0$ & $-35.4$ \\
$002$ & $0$ & $1$ & $-14.4$ \\
$020$ & $0$ & $1$ & $-13.1$ \\
$200$ & $0$ & $1$ & $-15.5$ \\
$021$ & $1$ & $1$ & $-25.6$ \\
$120$ & $1$ & $1$ & $-22.6$ \\
$121$ & $2$ & $1$ & $-35.1$ \\
$201$ & $1$ & $1$ & $-28.0$ \\
$210$ & $1$ & $1$ & $-26.0$ \\
$211$ & $2$ & $1$ & $-41.2$ \\
$012$ & $1$ & $1$ & $-24.9$ \\
$102$ & $1$ & $1$ & $-23.9$ \\
$102$ & $1$ & $1$ & $-23.9$ \\
$112$ & $2$ & $1$ & $-37.3$\\
$022$ & $0$ & $2$ & $-27.5$\\
$202$ & $0$ & $2$ & $-29.9$\\
$220$ & $0$ & $2$ & $-28.6$\\
$222$ & $0$ & $3$ & $-43.0$\\
$222$ & $0$ & $3$ & $-43.0$\\
$221$ & $1$ & $2$ & $-41.1$\\
$212$ & $1$ & $2$ & $-40.4$\\
$122$ & $1$ & $2$ & $-37.0$\\
\noalign{\smallskip}\hline
\end{tabular}
\end{table}


\begin{table}
\caption{Parameters for computing the vector field for the genetic toggle problem \cite{Aurell}.}
\label{Table:param}     
\centering
\begin{tabular}{llllllllll}
\hline\noalign{\smallskip}
$RT$ & $S_{CI}$ & $S_{Cro}$ & $R_{RM}$ & $R_{RM}^{u}$ & $R_R$ & $\tau_{CI}$ & $\tau_{Cro}$ & $\Delta G_{CI}$ & $\Delta G_{Cro}$\\
\noalign{\smallskip}\hline\noalign{\smallskip}
$0.617$ &  $1$  & $20$ & $0.115$ & $0.01045$ & $0.30$ & $2943$ & $5194$ & $-11.1$ & $-7.0$\\
\noalign{\smallskip}\hline
\end{tabular}
\end{table}




\end{document}